\documentclass{amsart}

\usepackage{amsmath}
\usepackage{amssymb}
\usepackage{amsthm}
\usepackage{color}
\usepackage{graphicx}
\usepackage{psfrag}
\usepackage{lipsum}% http://ctan.org/pkg/lipsum

\makeatletter
\g@addto@macro{\endabstract}{\@setabstract}
\newcommand{\authorfootnotes}{\renewcommand\thefootnote{\@fnsymbol\c@footnote}}%
\makeatother

\newcommand{\geqs}{\geqslant}
\newcommand{\leqs}{\leqslant}

\newcommand{\nhat}{\mathbf{n}}
\newcommand{\Rb}{{\mathbb R}}
\newcommand{\Z}{\mathbb{Z}}

\DeclareMathOperator{\curl}{curl}

\newtheorem{theorem}{Theorem}[section]

\newtheorem{prop}[theorem]{Proposition}

\theoremstyle{definition}
\theoremstyle{remark}

\numberwithin{equation}{section}

\begin{document}
\begin{center}
  {\sc
  Bifurcations analysis of the twist-Fr\'eedericksz \\ transition in a nematic liquid-crystal cell with \\ pre-twist
  boundary conditions: the asymmetric case}\footnotemark \par \bigskip

\footnotetext{F.P. da C. and J.T.P. were partially supported by FCT/Portugal through
UID/MAT/04459/2013.}

  \normalsize
  \authorfootnotes
  F.P. da Costa%\footnote{}
  \textsuperscript{2}, 
  M.I. M\'endez\footnote{Portions of this paper first appeared in the dissertation submited by the second author 
   to the \emph{Universidad Nacional de Educaci\'on a Distancia,\/} 
  Madrid, Spain, in June 2014, as part of the requirements for the  degree of
  \emph{M\'aster en Matem\'aticas Avanzadas.\/}}\textsuperscript{3},
  J.T. Pinto%\footnote{Author A}
  \textsuperscript{4} \par \bigskip

  \textsuperscript{2}Departamento de Ci\^encias e Tecnologia, Universidade Aberta, Lisboa, Portugal, and
Centro de An\'alise Matem\'atica, Geometria e Sistemas Din\^amicos, Instituto Superior T\'ecnico,
Universidade de Lisboa, Lisboa, Portugal 
\email{fcosta@uab.pt} \par
 \textsuperscript{3}Departamento de Matem\'aticas, IES Antonio L\'opez, Getafe (Madrid), Spain
\email{misabel.mendezaller@educa.madrid.org} \par 
  \textsuperscript{4}Departamento de Matem\'atica and 
Centro de An\'alise Matem\'atica Geometria e Sistemas Din\^amicos, 
Instituto Superior T\'ecnico, Universidade de Lisboa, Lisboa, Portugal
\email{jpinto@math.tecnico.ulisboa.pt} \par \bigskip

\end{center}

%%%%%%%%%%%%%%%%%%%%%%%%%%%%%%%%%%%%%%%%%%%%%%%%%%%%%%%%%%%%%%%%%%%%%%%%%%%%%%%%%%%%%%%%%%%%%%%%%%%%%%%%%%%%%%%%%%%%

%\begin{document}

%    General info
\begin{center}
\noindent
{\small December 10, 2015}
\end{center}

\subjclass[2010]{Primary 34B15, 34C23; Secondary  78A55}
\keywords{Two-points boundary value problems, Bifurcations, Time-maps, Fr\'eedericksz transition, Nematic liquid crystal cell}

%=============================================================================
%                                                                                 Abstract
%=============================================================================

\begin{abstract}
In the paper [Eur.\ J.\ of\ Appl.\ Math. \textbf{20}, (2009)  269--287] by da Costa et al. the twist-Fr\'eedericksz 
transition in a nematic liquid crystal one-dimensional cell of lenght $L$ was studied imposing an antisymmetric net 
twist Dirichlet condition at the cell boundaries. In the present paper we extend that study to the more general case of 
net twist Dirichlet conditions without any kind of symmetry restrictions. 
We use phase-plane analysis tools and appropriately defined time-maps 
to obtain the bifurcation diagrams of the model when $L$ is the bifurcation parameter, and related these diagrams with the
one in the symmetric situation. The stability of the bifurcating solutions is investigated by applying
the method of Kenjiro Maginu [J.\ Math.\ Anal.\ Appl. \textbf{63}, (1978) 224--243].
\end{abstract}

%\maketitle

%==================================================================================
%                                          Introduction
%==================================================================================

\section{Introduction}

In the operation of liquid crystal devices the phenomena of Fr\'eedericksz transitions in nematic liquid crystal cells are
of paramount technological importance \cite[Chapter 5]{yw} and give rise to interesting and challenging mathematical problems 
\cite[Section 3.4]{s}.

A nematic liquid crystal cell is basically a thin layer (a few microns) of a nematic liquid crystal contained
between two glass plates whose inner surface is chemically treated in such a way as to force a certain allignment 
(anchoring) of the 
rod like nematic liquid crystal molecules lying  close to the cell boundary. This surface alignment induces an allignment in 
the liquid crystal molecules filling the cell bulk so that the total free energy is minimized. 

When an exterior electric or magnetic field is applied to the cell a competition takes place
between the reorienting effects of the field and the allignment imposed by the surface anchoring.
Minimization of the total free energy (field, elastic bulk, and anchoring) then forces 
a reallignment of the molecules in the cell bulk (and, in the case of the so called weak anchoring conditions, 
also of those at the cell surface \cite{cgmp}) when the field intensity increases above a threshold value
dependent  on the physico-chemical characteristics of the device. This bifurcation phenomenon
 is called Fr\'eedericksz transition, in honor of the soviet physicist who discovered it \cite{fz}.

If we model the rod-like nematic liquid crystal molecules by a ``director vector field'' $\nhat,$
with $\|\nhat\|\equiv 1,$
a system with strong anchoring of the molecules at the cell surface occupying a region $\Omega$ 
has a total free energy of the director field given by
\[
\int_{\Omega} w( \nhat, \nabla\nhat ) ,
\]
where the free-energy density $w$ embodies the
competition between the energy cost of distortions of the director
field versus the energy reduction associated with aligning parallel (or perpendicular) to
the magnetic field, and is givem by
\[
2 w =
K_1 \left( \operatorname{div} \nhat \right)^2 +
K_2 \left( \nhat \cdot \curl \nhat \right)^2 +
K_3 \left\| \nhat \times \curl \nhat \right\|^2 -
\mu_0 \Delta \chi ( \mathbf{H} \cdot \nhat )^2 ,
\]
where $K_1$, $K_2$, and $K_3$ are phenomenological elastic constants, $\mu_0$ is the
free-space magnetic permeability, $\Delta\chi =
\chi_{\scriptscriptstyle\parallel} - \chi_{\scriptscriptstyle\perp}$
is the difference between the diamagnetic susceptibilities parallel to
versus perpendicular to the director, and $\mathbf{H}$ is the
(constant) applied magnetic field. See, e.g., \cite{s}.
 
\begin{figure}[h]\label{figLCCell}
\begin{center}
\psfrag{d}{{\Large $d$}}
\psfrag{H}{{\Large $\mathbf{H}$}}
\psfrag{n}{{\Large $\mathbf{n}$}}
\psfrag{phi0}{$-\phi_0$}
\psfrag{phi1}{$\phi_1$}
%\psfrag{S1}{$\Sigma_1(t)$}
\includegraphics[scale=0.35]{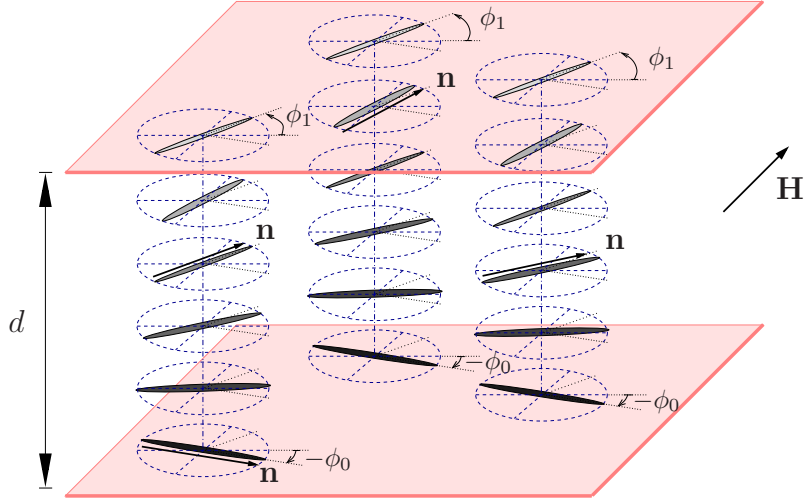}
\end{center}
\caption{Geometry of the liquid crystal cell with asymmetric pre-twist. The director $\mathbf{n}$ orientation inside the cell
corresponds to the situation in branch $C_r$ (with $k=0$) in Section~3.1 below.}
\end{figure}

We consider the geometry of the twist-Fr\'eedericksz transition, with
an asymmetric pre-twist at the boundary. Thus we consider a thin slab of nematic
liquid crystal bounded by two
parallel planes a distance $d$ apart from each other, unbounded and
extending to infinity in any direction parallel to these planes.
Define a positively oriented orthogonal coordinate system $(x, y, z)$
such that $z$ is perpendicular to the bounding planes.  Let the
director field be represented by
\begin{equation}
  \label{eqn:nhat}
  \mathbf{n} = ( \cos\phi(z,\tau), \sin\phi(z, \tau), 0 ) ,
\end{equation}
where $\phi$ denotes the (twist) angle of the director. We will assume
that in the liquid crystal cell the director is fixed in opposing
orientations $-\phi_0$ and $\phi_1$ at the two opposing planes
bounding the device in the $z$ direction. This induces a net twist of
the director vector field $\mathbf{n}$ across the cell (see
Figure~\ref{figLCCell}).

%
%%%%%%%%%%%%%%%%%%%%%%%%%%%%%%%%%%%%%%%%%%%%%%%%%%%%%%%%%%%%%%%%%%%%%%%%
%
%
%
%
%
%
%\begin{figure}[h]\label{figLCCell}
%
%
%
%
%
%
%%%%%%%%%%%%%%%%%%%%%%%%%%%%%%%%%%%%%%%%%%%%%%%%%%%%%%%%%%%%%%%%%%%%%%%%
%
We will consider a magnetic vector field $\mathbf{H}$ applied along
the constant direction $(0, 1, 0)$ with intensity
$H=\|\mathbf{H}\|$ and are interested in studying the effect it
induces in the stationary director distribution, according to the
Ericksen-Leslie theory \cite{s}.

In terms of the angle representation \ref{eqn:nhat}, the simplest model
for the dynamics of the director field in the absence of
flow is the gradient flow on the \emph{free energy} of
the system and, in dimensionless form, the initial-boundary value problem governing
the behaviour of the director field is then
\begin{eqnarray}
  & & \frac{\partial \phi}{\partial s} = \frac{\partial^2
    \phi}{\partial \zeta^2} + \lambda \sin\phi \cos\phi, \qquad \quad
  (s,\zeta) \in \Rb^+ \times (0,1) \label{eq1+} \\ 
  & & \phi(\cdot,0) \, = \, - \phi_0 , \qquad \phi(\cdot,1) \, = \,
  \phi_1 , \label{eq2+}\\
  & & \phi(0,\cdot) \, = \, \phi_{\rm initial} \label{eq1##}
\end{eqnarray}
where
\begin{equation}
s := \frac{K_2}{\gamma_1d^2} \tau\;\;, \quad \zeta:=\frac{z}{d}\;\;, \quad
\lambda := \frac{\mu_0\Delta\chi H^2d^2}{K_2}, \label{adimensional}
\end{equation}
with all the material parameters are positive for our system of interest.
Observe that the dimensionless control parameter $\lambda$ is
proportional to the square of the magnetic field strength.

The associated equilibrium problem is given by
\begin{eqnarray}
  \frac{d^2\phi}{d\zeta^2} + \lambda \sin\phi \cos\phi &=& 0 , \qquad
  0 < \zeta < 1  \label{*} \\
  \phi(0) = - \phi_0 , \quad \phi(1) = \phi_1. & & \label{**}
\end{eqnarray}

In the classical twist-Fr\'eedericksz-transition problem, we have
$\phi_0=\phi_1=0$ the system possesses a simple symmetry,
$\phi(\zeta) \leftrightarrow - \phi(\zeta),$ and the ground-state
solution ($\phi=0$, which is invariant under this symmetry) loses
stability to a pair of symmetric solutions at a
pitchfork bifurcation at $\lambda_c = \pi^2.$ 

In \cite{CGGP} a system with antisymmetric pre-twist 
($\phi_0=\phi_1\neq 0$)  was studied. We no longer have the simple symmetry above.  
The problem still possesses $\Z_2$ symmetry,
however it is now of the form $\phi(\zeta) \leftrightarrow - \phi(1-\zeta)$.
The ground-state solution (which is invariant under this symmetry) is no
longer uniform. The problem still has a classical pitchfork bifurcation diagram,
with the symmetric solution branch bifurcating at a value
$\lambda_c$, which is necessarily greater than $\pi^2,$ as was showed in \cite{CGGP}.
Observe that the antisymmetric nature of the boundary data is crucial to this scenario.

In the present paper we consider the asymmetric case ($\phi_0\neq \phi_1,$ both nonzero).
We conclude that no pitchfork bifurcation points remain: the pitchforks that had not been broken 
in the passage from the classical twist cenario to the 
antisymmetric one, are now broken when the $\phi_0$ becomes different from $\phi_1,$ and the result
is a bifurcation diagram with only saddle-node bifurcation points, branches emanating from them, and 
single nonbifurcating branch of solutions. 

The approach will be based on the time maps and phase-plane methods developed in
\cite{CGGP} for the antisymmetric case.
The stability of these branches is also studied by applying the results of \cite{maginu78}, also based
on the behaviour of time-maps, which allows the classification of the stationary solution branches
as stable, asymptotically stable, or unstable. A more detailed study of the stability indices of the equilibria
and the characterization of their connecting orbits will be the subject of a future paper.

%==================================================================================
%                                          Preliminaries
%==================================================================================

\section{Preliminaries}

We will be concerned with the stationary
solutions to (\ref{eq1+})-(\ref{eq1##}), i.e., solutions of (\ref{*})-(\ref{**}). Consider the change
of variables $t = t(\zeta) := \sqrt{\frac{\lambda}{2}} \left( \zeta -
  \frac{1}{2} \right),$ and let $\zeta(t)$ be its inverse function.
Let 
\begin{equation}
L := \sqrt{\frac{\lambda}{8}}. \label{L}
\end{equation}
Then, $\phi(\zeta)$
is a solution of (\ref{*})-(\ref{**}) iff $x(t) :=
\phi(\zeta(t))$ is a solution of
\begin{eqnarray}
  & & \left\{\begin{array}{l}x' = y \\
      y' = -\sin 2x \end{array} \right. \label{eq1} \\ 
  & & \nonumber \\
  & & x(-L) = - \phi_0, \quad x(L)=\phi_1, \label{eq2}
\end{eqnarray}
where  $\phi_0, \phi_1 \in (0, \frac{\pi}{2}),$ and 
$(t,x,y)\in [-L,L]\times[-\pi/2, \pi/2]\times\Rb.$
The bifurcation parameter is now $L>0$.  Note that $L\propto H.$ 
We shall treat the independent variable $t$ in (\ref{eq1})-(\ref{eq2}) as the ``time''
of the dynamical system associated with (\ref{eq1+}). Note that this ``time'' corresponds
to the original spacial variable $\zeta$ and not to the original time $s$.

The study of the bifurcation structure of solutions to (\ref{eq1})--(\ref{eq2}) when $\phi_0=\phi_1$ was done
in \cite{CGGP}. We now consider the general case, where no relation between 
the values of $\phi_0$ and $\phi_1$ is imposed. As in \cite{CGGP}, we shall use the tools of time-maps
and phase-plane analysis.

\begin{figure}[h]
\begin{center}
\psfrag{x}{$x$}
\psfrag{y}{$y$}
\psfrag{-pi}{\scriptsize{$\left(-\frac{\pi}{2},0\right)$}}
\psfrag{pi}{\scriptsize{$\left(\frac{\pi}{2},0\right)$}}
\psfrag{x=+}{$x=\phi_1$}
\psfrag{x=-}{$x=-\phi_0$}
\includegraphics[scale=0.41]{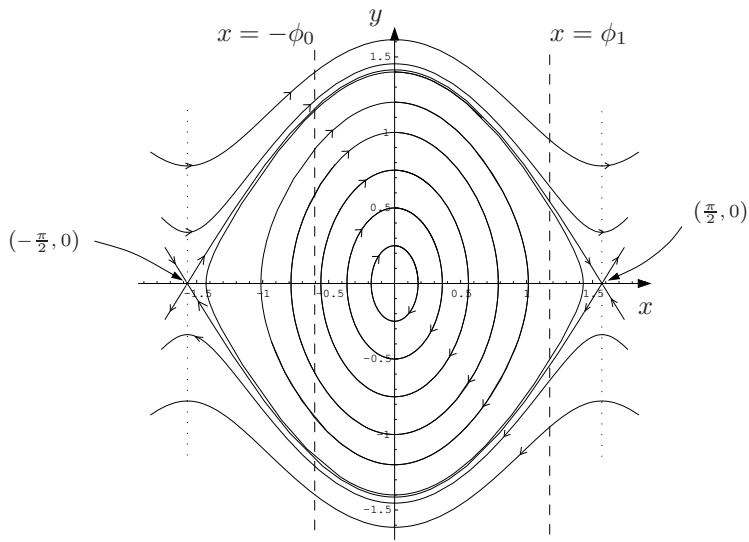} 
\end{center}
\caption{Orbits of \eqref{eq1} with the Dirichlet boundary condition \eqref{eq2}
marked by the dashed vertical lines (with $0<\phi_0<\phi_1<\pi/2.$)} \label{figOrbits}
\end{figure}

The phase portrait of \eqref{eq1} is presented in Figure~\ref{figOrbits}.

%%%%%%%%%%%%%%%%%%%%%%%%%%%%%%%%%%%%%%%%%%%%%%%%%%%%%%%%%%%%%%%%%%%%%%%%
%
%
%
%
%
%
%\begin{figure}[h]
%
%
%
%
%
%
%%%%%%%%%%%%%%%%%%%%%%%%%%%%%%%%%%%%%%%%%%%%%%%%%%%%%%%%%%%%%%%%%%%%%%%%

For studying solutions to \eqref{eq1}-\eqref{eq2} we need to consider some time-maps measuring the time spent by 
the orbits. These maps are easily obtained from the fact that \eqref{eq1}-\eqref{eq2} is a conservative system with
energy 
\begin{equation}
 V(x,y) = y^2 -\cos 2x, \label{energyV}
\end{equation}
which means that its orbits are subsets of the level sets of this function.

Let $\alpha\in\left(0, \frac{\pi}{2}\right)$ and denote by $\gamma_{\alpha}$ 
the orbit that, at time $t=0$, intersects the $x$-axis at $(\alpha, 0).$ Clearly $\gamma_{\alpha}$ is
a periodic orbit (cf. Figure~\ref{figOrbits}). Let $P(\alpha)$ be its period and define
\begin{equation}
  T(\alpha) := \frac{1}{4}P(\alpha) = \int_0^{\alpha}
  \frac{dx}{\sqrt{\cos 2x-\cos 2\alpha}}, \label{T} 
\end{equation}
where the second equality arises from the symmetry of the system with respect to reflexions in the $x-$ and $y-$axis.
Thus, $T(\alpha)$ is the time it takes for the point of intersection of $\gamma_{\alpha}$ with the $y-$axis (which, by conservation
of $V$ along orbits, we easily conclude to be the point $(0, \beta)$ with $\beta=\sqrt{2}\sin \alpha$) to reach the $x-$axis (at the point
$(\alpha, 0),$ by construction).

As in \cite{CGGP}, two other time maps will be needed. The time-map 
\begin{equation}
  T_1(\alpha,\phi) := \int_0^{\phi}
  \frac{dx}{\sqrt{\cos 2x-\cos 2\alpha}}, \label{T1} 
\end{equation}
that measures the time spent by the point of intersection of the orbit $\gamma_{\alpha}$ with the $y-$axis to reach the line $x=\phi\leqs \alpha.$
Clearly $T_1(\phi, \phi) = T(\phi).$
We will also consider a map $T_2$ analogous to $T_1$ but relevant 
for orbits crossing the $y-$axis on or above the heteroclinic orbit $\gamma_h$ connecting $(-\pi/2,0)$ to $(\pi/2, 0)$, 
i.e., at a point $(0, \beta)$
with $\beta \geqs \sqrt{2}$, namely
\begin{equation}
  T_2(\beta,\phi) := \int_0^{\phi}
  \frac{dx}{\sqrt{\beta^2+\cos 2x-1}}. \label{T2} 
\end{equation}
We can continuously extend this map to values $\beta<\sqrt{2}$ by $T_2(\beta, \phi) := T_1(\alpha(\beta),\phi)$, where $\alpha(\beta)$ is defined to be the unique
value of $\alpha$ for which the points $(\alpha,0)$ and $(0,\beta$) are on the same orbit. Since the orbits are contained in the level sets of $V$,
a brief inspection of Figure~\ref{figOrbits} allow us to conclude that
 $\beta\mapsto \alpha(\beta)$ is a monotonically increasing function and thus, for each fixed $\phi$, there is a smaller $\beta$ for which $T_2(\beta, \phi)$ is
 defined, which is the value $\beta_{\phi}$ for which $(\beta_{\phi}, 0)$ and $(\phi, 0)$ are on the same orbit. For $\beta$ below $\beta_{\phi}$
 no orbit  satisfies the boundary condition at $t=L$. 

 Our analysis depend heavily on the following monotonicity properties of the time-maps defined above. 
 A proof of these results can be checked in \cite{CGGP}.

\begin{prop} \label{monotonicity}
Let $\alpha\in\left(0, \frac{\pi}{2}\right),$ $\phi\in(0,
  \alpha),$ and $\beta\geqslant\beta_{\phi}.$ Then, 
\begin{description}
 \item[1] the time-map $T: \left(0, \frac{\pi}{2}\right)\to (0, +\infty)$
  defined by (\ref{T}) is strictly increasing and converges to
  $\frac{\pi}{2\sqrt{2}}$ as $\alpha\to 0$, and to $+\infty$ as $\alpha \to \frac{\pi}{2}.$
  
 \item[2] for each fixed $\phi$ the time-maps $T_1(\cdot,\phi)$ and $T_2(\cdot,\phi)$, defined by
 (\ref{T1}) and (\ref{T2}), respectively, are  strictly decreasing. 
The same holds true for $T_2(\cdot, \frac{\pi}{2}).$
\end{description}
\end{prop}

%==================================================================================
%                                          Phase-plane analysis
%==================================================================================

\section{Bifurcation analysis}

The study of \eqref{eq1}-\eqref{eq2} in the symmetric case $\phi_0=\phi_1$ was done in \cite{CGGP} and will serve as a 
guide to our present study. In the symmetric case a special role is played by the solutions of \eqref{eq1}-\eqref{eq2}
that additionally satisfy the homogeneous Neumann boundary condition $y(-L)=y(L)=0$ (note that $y=x'$). The values
of $L$ for which these solutions occur were termed ``critical'' (cf. \cite[Figure 4 and Table 1]{CGGP}) and 
are pitchfork bifurcation points of the system \cite[Figure 9]{CGGP}. The orbits corresponding to these
values of $L$ were denoted by $\gamma^*.$

Due to the symmetry of the vector field of \eqref{eq1} and the asymmetry of the boundary condition \eqref{eq2}
there are no solutions to \eqref{eq1}-\eqref{eq2} satisfying homogeneous Neumann boundary conditions at \emph{both}
$t=-L$ and $t=L.$ However, there are solutions that satisfy such a condition at one, or the other, of the end points of
the time interval. Although these do not correspond to bifurcation solutions, and the corresponding values of $L$
are not bifurcation points of \eqref{eq1}-\eqref{eq2}, they are important solutions that help us 
to organize the information and construct the bifurcation diagram in the asymmetric case, and relate it
with the symmetric case already studied.

The two asymmetric cases $\phi_0<\phi_1$ and $\phi_0>\phi_1$ give rise to different bifurcation diagrams and will
be studied separately below. Since the approach for both cases is the same, we will present the first one 
in a more detailed way, and for the second will just briefly refer to the corresponding results.

\subsection{Case $\phi_0<\phi_1$}

 \subsubsection{The ``critical'' cases}\label{311}

\begin{figure}[h]
\begin{center}
\psfrag{a}{\scriptsize{$(a)$}}
\psfrag{b}{\scriptsize{$(b)$}}
\psfrag{x}{$x$}
\psfrag{y}{$y$}
\psfrag{-p}{\scriptsize{$x=-\phi_0$}}
\psfrag{p}{\scriptsize{$x=\phi_1$}}
\psfrag{-pi}{$-\frac{\pi}{2}$}
\psfrag{pi}{$\frac{\pi}{2}$}
\psfrag{ga}{$\gamma_*$}
\psfrag{gb}{$\gamma^*$}
\includegraphics[scale=0.23]{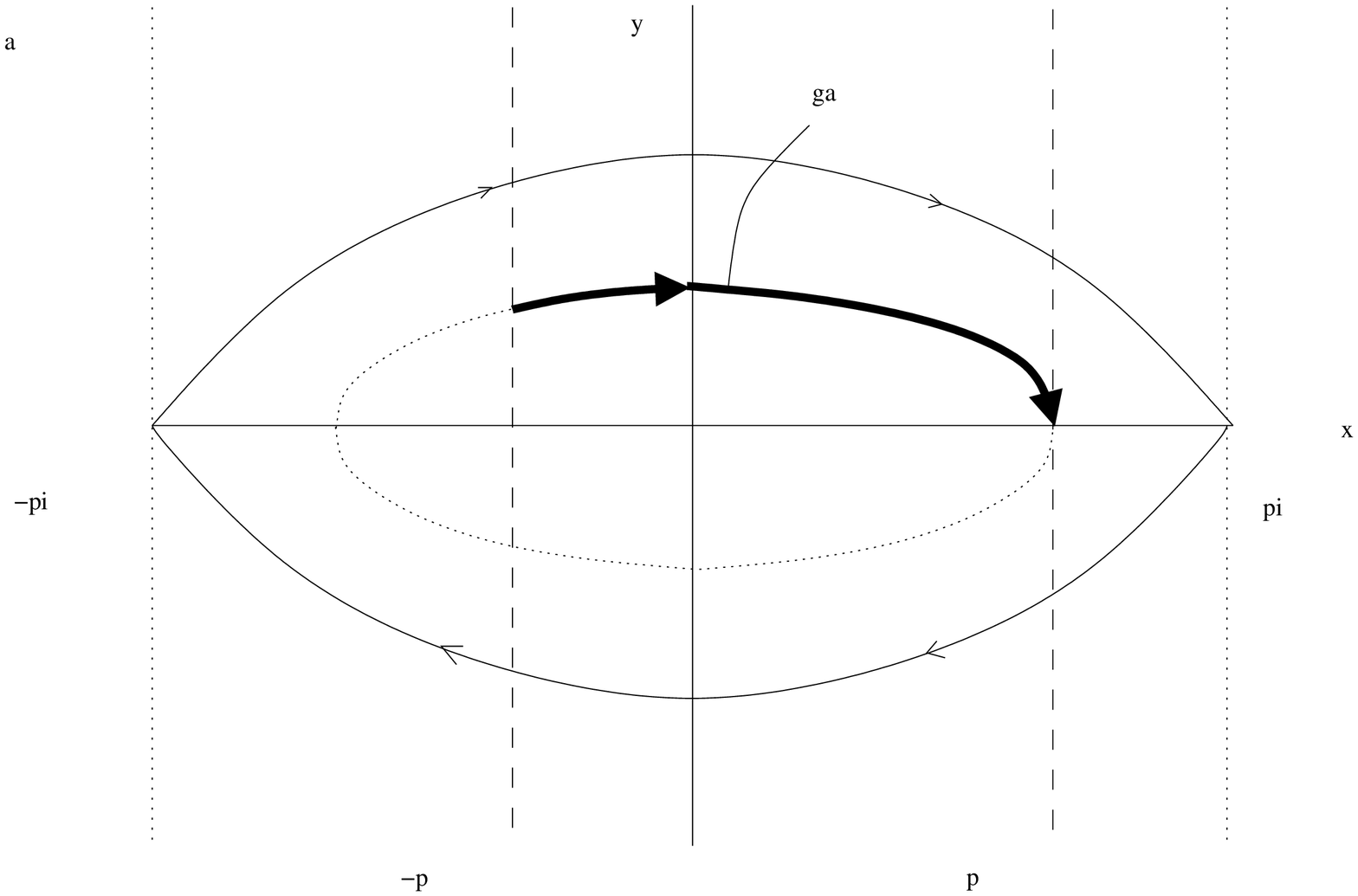}\qquad\includegraphics[scale=0.23]{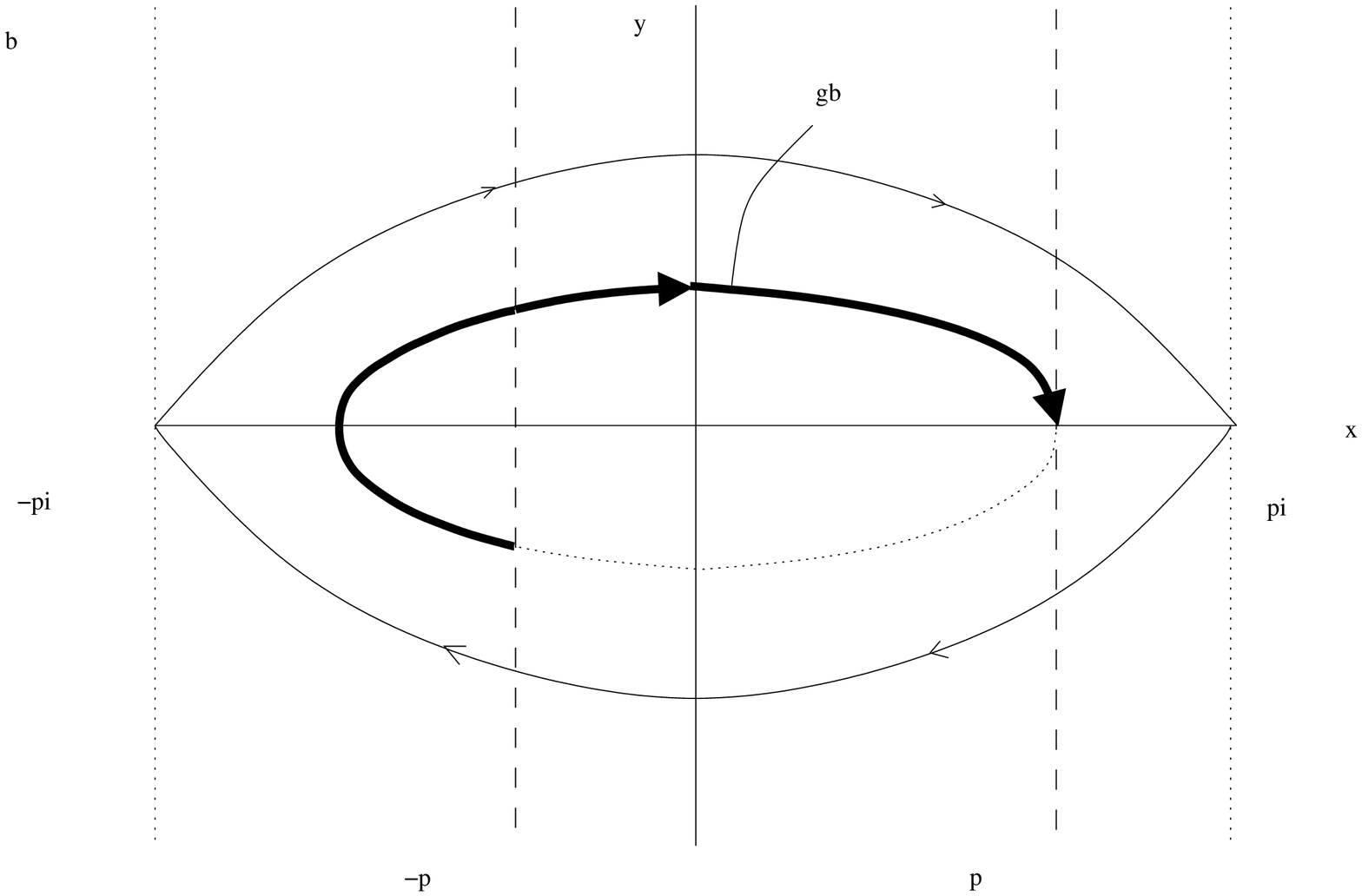}
\end{center}
\caption{Two ``critical'' orbits of (\ref{eq1})-(\ref{eq2}) with $\phi_0<\phi_1$: (a) $\gamma_*$ 
 when $2L=T_*:=T(\phi_1)+T_1(\phi_1,\phi_0)$;
 (b) $\gamma^*$ 
 when $2L=T^*:=3T(\phi_1)-T_1(\phi_1,\phi_0)$} \label{figcritic}
\end{figure}

 Let $\gamma_*$ be the orbit of \eqref{eq1}-\eqref{eq2} that satisfies the additional homogeneous Neumann
 condition $y(L)=0.$ See Figure~\ref{figcritic}(a). It is clear from this figure and from the definition
 of the time-maps in the previous section that the time spent by $\gamma_*$ is $T_*:=T(\phi_1)+T_1(\phi_1,\phi_0).$ 
 Since the total time spent by every orbit is $2L$, the corresponding half-length $L$ is $L_*=\frac{1}{2}T_*.$
 
 In a similar way, the orbit that satisfies the homogeneous Neumann boundary condition $y(-L)=0$
 will be denoted by $\gamma^*.$ See Figure~\ref{figcritic}(b). The time spent in by this orbit is 
 $T^*:=3T(\phi_1)-T_1(\phi_1,\phi_0)$, and the corresponding half-lenght of the interval is $L^*=\frac{1}{2}T^*.$

%%%%%%%%%%%%%%%%%%%%%%%%%%%%%%%%%%%%%%%%%%%%%%%%%%%%%%%%%%%%%%%%%%%%%%%%
%
%
%
%
%
%
%
%
%
%
%
%
%%%%%%%%%%%%%%%%%%%%%%%%%%%%%%%%%%%%%%%%%%%%%%%%%%%%%%%%%%%%%%%%%%%%%%%%

From the definitions of the time-maps we easily observe that, 
\begin{equation}
 T_* =T(\phi_1) +T_1(\phi_1,\phi_0) < 2T(\phi_1) < 3T(\phi_1) - T_1(\phi_1, \phi_0) = T^*, \label{desigualdadeT}
\end{equation}
and these the inequalities turn to equalities if $\phi_0=\phi_1$, which, as already pointed out, 
was the case considered in \cite{CGGP}.

 By analogy to the terminology used in \cite{CGGP} for the symmetric case we shall call these solutions, orbits, etc.,
 ``critical'', although, as we shall see, they do not correspond to any critical feature in the bifurcation diagrams.
 However, they will be very useful for the remaining constructions. In particular, as a matter of terminology
 and when appropriate, we will keep
 denoting by subcritical [resp., supercritical]
 those situations with values of $L$ smaller [resp., larger] than $L^*$ ou $L_*$.

\subsubsection{The subcritical case relative to $\gamma_*$}\label{312}

By Proposition~\ref{monotonicity} it is clear that the function 
$\left(\phi_1, \frac{\pi}{2}\right)\ni\alpha\mapsto T_A(\alpha):= T_1(\alpha, \phi_1)
+T_1(\alpha, \phi_0)$ is monotonically decreasing and $T_A(\alpha)\uparrow T_*$ as 
$\alpha\downarrow \phi_1$. The corresponding orbit of \eqref{eq1}-\eqref{eq2} is a subset of the level 
set $V(\alpha,0)$ of $V.$ Since the time it spents is $2L=T_A(\alpha)< T_*,$ we call it subcritical 
relative to $\gamma_*$. Using the relation between 
the time-maps $T_1$ and $T_2$ we can extend this approach to orbits intersecting the $y-$axis above 
the heteroclinic orbit $\gamma_h.$ The time taken by these
orbits is also smaller than $T_*$ and decreases as the ordinate of the intersection point increases. 

\begin{figure}[h]
\begin{center}
\psfrag{a}{\scriptsize{$(a)$}}
\psfrag{b}{\scriptsize{$(b)$}}
\psfrag{x}{$x$}
\psfrag{y}{$y$}
\psfrag{-p}{\scriptsize{$x=-\phi_0$}}
\psfrag{p}{\scriptsize{$x=\phi_1$}}
\psfrag{-pi}{$-\frac{\pi}{2}$}
\psfrag{pi}{$\frac{\pi}{2}$}
\psfrag{ga}{$\gamma_*$}
\psfrag{gb}{$\gamma^*$}
\includegraphics[scale=0.23]{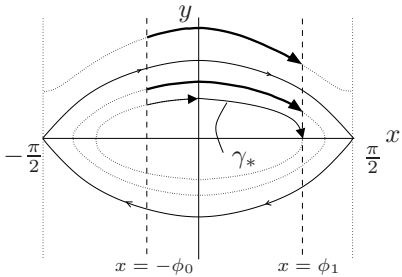}
\end{center}
\caption{Two  orbits of (\ref{eq1})-(\ref{eq2}), with $\phi_0<\phi_1$, subcritical relative to the orbit
$\gamma_*$.} \label{figsubcritic}
\end{figure}

In Figure~\ref{figsubcritic} we present two of these orbits subcritical relative to $\gamma_*$, 
together with the critical orbit $\gamma_*.$
The monotonicity of the time-maps imply that, for each $L\in (0, \frac{1}{2}T_*)$ there is a 
single subcritical solution to \eqref{eq1}-\eqref{eq2}.

%%%%%%%%%%%%%%%%%%%%%%%%%%%%%%%%%%%%%%%%%%%%%%%%%%%%%%%%%%%%%%%%%%%%%%%%
%
%
%
%
%
%
%
%
%
%
%
%
%%%%%%%%%%%%%%%%%%%%%%%%%%%%%%%%%%%%%%%%%%%%%%%%%%%%%%%%%%%%%%%%%%%%%%%%

\subsubsection{The supercritical case relative to $\gamma_*$} \label{313}

\begin{figure}[h]
\begin{center}
\psfrag{a}{\scriptsize{$(a)$}}
\psfrag{b}{\scriptsize{$(b)$}}
\psfrag{x}{$x$}
\psfrag{y}{$y$}
\psfrag{-p}{\scriptsize{$x=-\phi_0$}}
\psfrag{p}{\scriptsize{$x=\phi_1$}}
\psfrag{-pi}{$-\frac{\pi}{2}$}
\psfrag{pi}{$\frac{\pi}{2}$}
\psfrag{ga}{$\gamma_*$}
\psfrag{gb}{$\gamma^*$}
\includegraphics[scale=0.23]{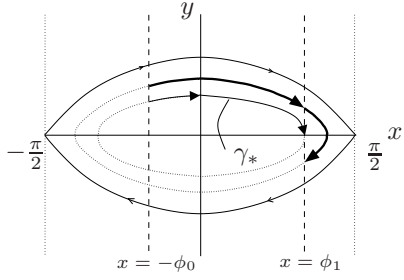}
\end{center}
\caption{An  orbit of (\ref{eq1})-(\ref{eq2}), with $\phi_0<\phi_1$, supercritical relative to the orbit
$\gamma_*$.} \label{figsupcritic}
\end{figure}

Consider again $\alpha\in \left(\phi_1, \frac{\pi}{2}\right)$ and the level set $V(\alpha, 0)$ of $V.$ 
For $\alpha>\phi_1$ but close to $\phi_1$, we take an orbit of (\ref{eq1})-(\ref{eq2}) close to $\gamma_*$
which have its end point with $y(L)<0$, as presented in Figure~\ref{figsupcritic}. A brief inspection of this figure
allow us to conclude that the time spent by this orbit is $T_{C_r}(\alpha) := 2T(\alpha) + T_1(\alpha,\phi_0) - T_1(\alpha, \phi_1)$
(the notation $T_{C_r}$ was used in \cite{CGGP} for a branch of solutions with a given symmetry relative to the origin.
We use the same notation here because our $C_r$ solutions will coincide with those of that paper when $\phi_0=\phi_1$; 
the same will be done for other solution branches further down the paper).  
Clearly $T_{C_r}(\alpha) \to T_*$ as $\alpha\to \phi_1.$

%%%%%%%%%%%%%%%%%%%%%%%%%%%%%%%%%%%%%%%%%%%%%%%%%%%%%%%%%%%%%%%%%%%%%%%%
%
%
%
%
%
%
%
%
%
%
%
%
%%%%%%%%%%%%%%%%%%%%%%%%%%%%%%%%%%%%%%%%%%%%%%%%%%%%%%%%%%%%%%%%%%%%%%%%

We shall prove that $T_{C_r}(\alpha) >T_*$ for  $\alpha >\phi_1$, thus providing a justification for calling these 
orbits supercritical
relative to $\gamma_*$. From the expression of $T_{C_r}$ and Proposition~\ref{monotonicity} we conclude that
\begin{equation}
 \frac{dT_{C_r}}{d\alpha}(\alpha) = T'(\alpha) + \frac{\partial T_1}{\partial \alpha}(\alpha, \phi_0) - \frac{\partial T_1}{\partial \alpha}(\alpha, \phi_1)
 > \frac{\partial T_1}{\partial \alpha}(\alpha, \phi_0) - \frac{\partial T_1}{\partial \alpha}(\alpha, \phi_1).\label{timebound}
\end{equation}
To conclude the sign of the right-hand side observe that
\[
% \begin{equation}
 \frac{\partial}{\partial \phi}\frac{\partial T_1}{\partial\alpha} = \frac{\partial}{\partial\alpha}\frac{\partial T_1}{\partial \phi} = 
 \frac{\partial}{\partial\alpha}\frac{1}{\sqrt{\cos 2\phi - \cos 2\alpha}} = - \frac{\sin 2\alpha}{(\cos 2\phi - \cos 2\alpha)^{3/2}} <0. \label{sign2derivative}
% \end{equation}
\]
From this inequality and the assumption that $\phi_0<\phi_1$ we infer that 
\[
% \begin{equation}
 \frac{\partial T_1}{\partial \alpha}(\alpha, \phi_0) > \frac{\partial T_1}{\partial \alpha}(\alpha, \phi_1),
% \end{equation}
\]
and plugging this into \eqref{timebound} gives that $T_{C_r}$ is strictly increasing with $\alpha$,  concluding the proof.

\subsubsection{The supercritical case relative to $\gamma^*$} \label{314}

\begin{figure}[h]
\begin{center}
\psfrag{a}{\scriptsize{$(a)$}}
\psfrag{b}{\scriptsize{$(b)$}}
\psfrag{x}{$x$}
\psfrag{y}{$y$}
\psfrag{O1}{$\Omega^-$}
\psfrag{O2}{$\Omega^+$}
\psfrag{-p}{\scriptsize{$x=-\phi_0$}}
\psfrag{p}{\scriptsize{$x=\phi_1$}}
\psfrag{-pi}{$-\frac{\pi}{2}$}
\psfrag{pi}{$\frac{\pi}{2}$}
\psfrag{ga}{$\gamma_*$}
\psfrag{gb}{$\gamma^*$}
\includegraphics[scale=0.23]{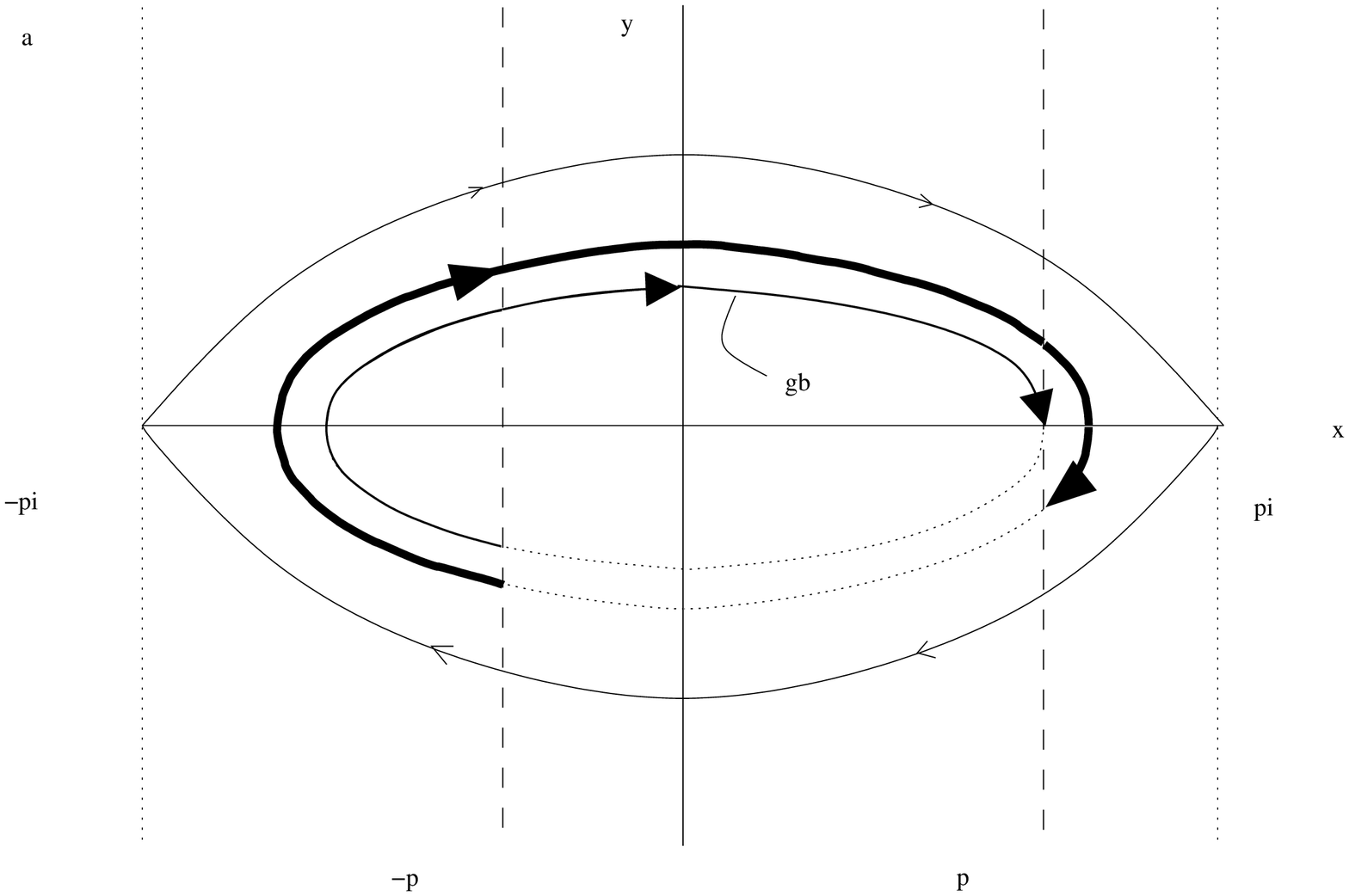} \quad \includegraphics[scale=0.23]{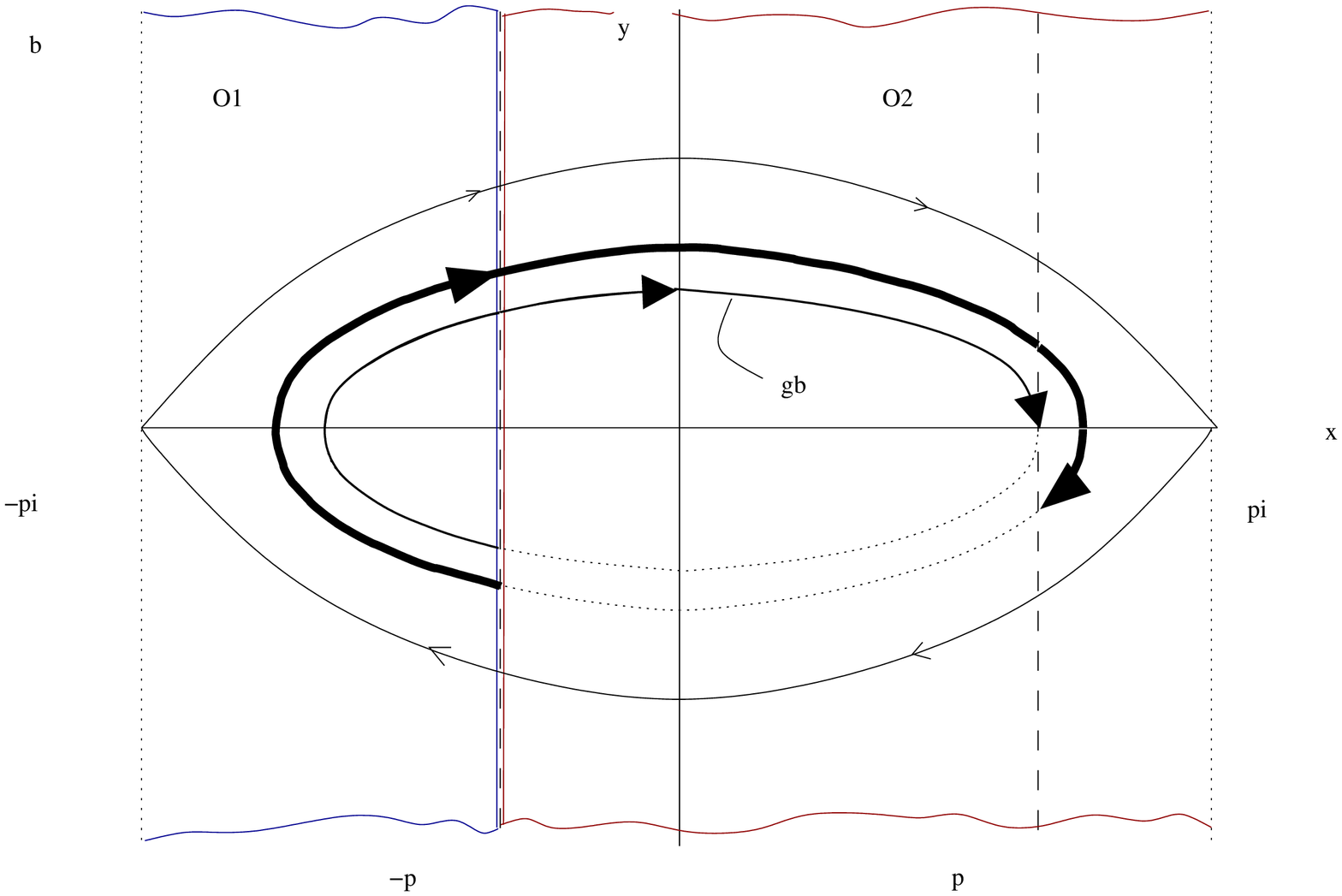}
\end{center}
\caption{(a) An  orbit of (\ref{eq1})-(\ref{eq2}), with $\phi_0<\phi_1$, supercritical relative to the orbit
$\gamma^*$. (b) The $\Omega^-$ and $\Omega^+$ regions.} \label{figsupcritic2}
\end{figure}

Consider an orbit in $V(\alpha, 0)$, with $\alpha\in \left(\phi_1, \frac{\pi}{2}\right),$ as represented in 
Figure~\ref{figsupcritic2}(a).
From this figure and the definition of the time-maps
we immediately conclude that the time spent to travel this orbit is 
$T_{D}(\alpha) := 4T(\alpha) - T_1(\alpha,\phi_0) - T_1(\alpha, \phi_1)$
(see subsection~\ref{313} for a justification of this notation).  

%%%%%%%%%%%%%%%%%%%%%%%%%%%%%%%%%%%%%%%%%%%%%%%%%%%%%%%%%%%%%%%%%%%%%%%%
%
%
%
%
%
%
%
%
%
%
%
%
%%%%%%%%%%%%%%%%%%%%%%%%%%%%%%%%%%%%%%%%%%%%%%%%%%%%%%%%%%%%%%%%%%%%%%%%

Clearly $T_{D}(\alpha) \to T_*$ as $\alpha\to \phi_1.$ We shall prove that $T_{D}(\alpha) > T^*.$
In order to prove this, consider the strips $\Omega^- := (-\pi/2, -\phi_0)\times\Rb,$
and $\Omega^+ := [-\phi_0, \pi/2)\times\Rb.$ Let ${\gamma^*}^{\pm} := \gamma^*\cap \Omega^{\pm}.$  
Denoting by $D$ an orbit of the type represented in Figure~\ref{figsupcritic2},
let also $D^{\pm} = D\cap \Omega^{\pm}.$ Since ${\gamma^*}^+ = \gamma_*,$ 
the time spent in ${\gamma^*}^+$ is equal to $T_*.$ 

Thus,  in $\Omega^+$ we just need to compare $T_*$ with the time spent by $D^+.$
But $D^+$ is really an orbit of type $C_r$ with $\alpha>\phi_1$ and thus, by the previous subsection,
$T_{D^+}(\alpha) >  T_*.$

In $\Omega^-$ we need to compare the time taken by the orbit $D^-$ with 
that taken by $\gamma^{*-},$ which a brief inspection to Figure~\ref{figsupcritic2}(b)
shows it is equal to $2 T(\phi_1) - 2T_1(\phi_1, \phi_0).$ Since
\[
 T_{D^-}(\alpha)  =  2T(\alpha) - 2T_1(\alpha, \phi_0),\quad \alpha \in (\phi_1, \pi/2),
\]
we have $\frac{\partial T_{D^-}}{\partial \alpha}(\alpha)  = 2T'(\alpha) - 2\frac{\partial T_1}{\partial \alpha}
(\alpha, \phi_0),$ and the monotonicity results in Proposition~\ref{monotonicity} imply that this derivative is positive, 
and thus $T_{D^-}(\alpha)>2 T(\phi_1) - 2T_1(\phi_1, \phi_0).$

Finally, from the above we have
\begin{eqnarray*}
T_{D}(\alpha) & = & T_{D^-}(\alpha)+T_{D^+}(\alpha) \\
              & > & 2 T(\phi_1) - 2T_1(\phi_1, \phi_0) + T_* = 3 T(\phi_1) - T_1(\phi_1, \phi_0)\\
              & = & T^*,
\end{eqnarray*}
as we wanted to prove.

\subsubsection{The subcritical case relative to $\gamma^*$} \label{315}

\begin{figure}[h]
\begin{center}
\psfrag{a}{\scriptsize{$(a)$}}
\psfrag{b}{\scriptsize{$(b)$}}
\psfrag{x}{$x$}
\psfrag{y}{$y$}
\psfrag{-p}{\scriptsize{$x=-\phi_0$}}
\psfrag{p}{\scriptsize{$x=\phi_1$}}
\psfrag{-pi}{$-\frac{\pi}{2}$}
\psfrag{pi}{$\frac{\pi}{2}$}
\psfrag{ga}{$\gamma_*$}
\psfrag{gb}{$\gamma^*$}
\includegraphics[scale=0.23]{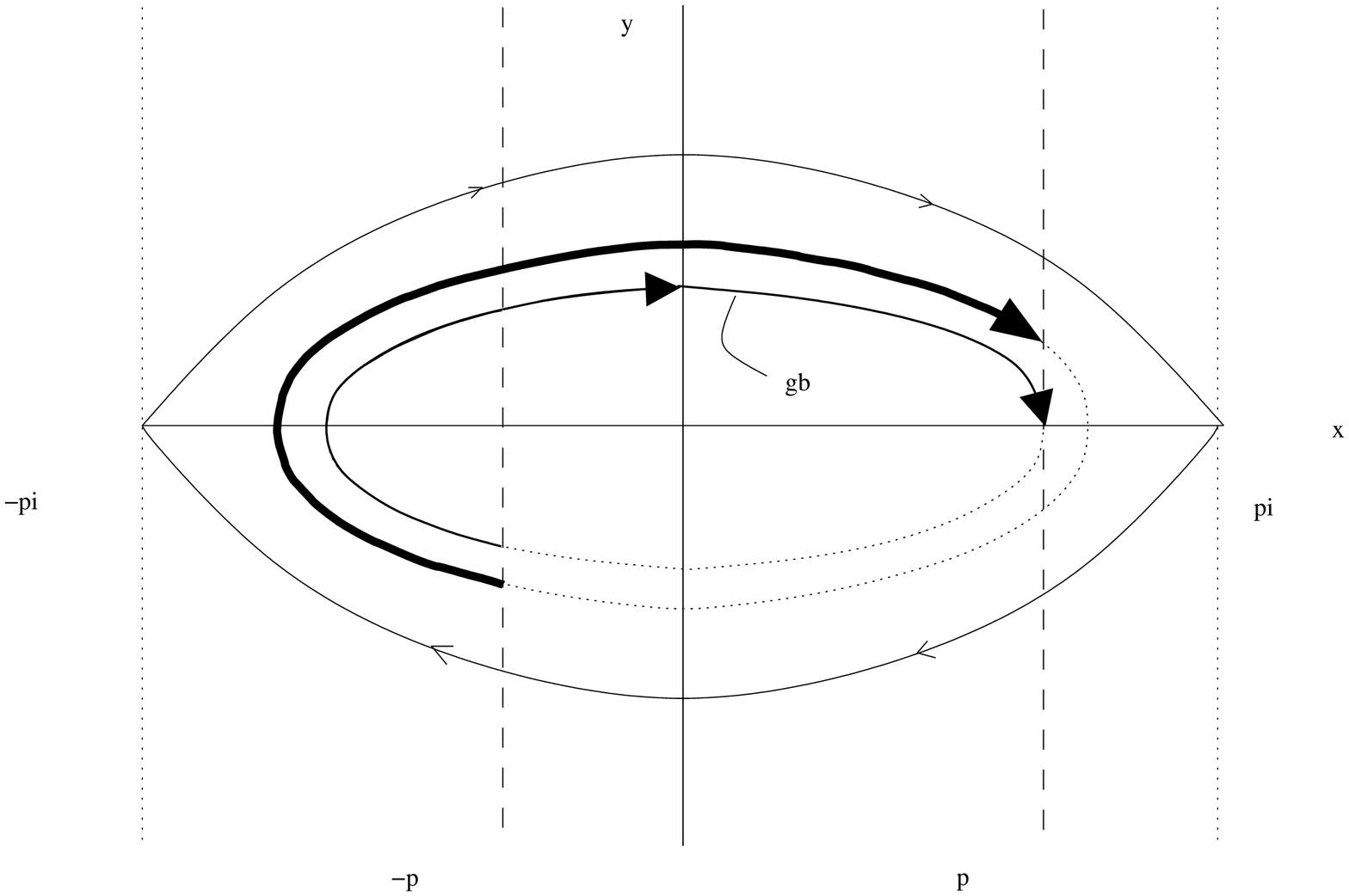}
\end{center}
\caption{An  orbit of (\ref{eq1})-(\ref{eq2}), with $\phi_0<\phi_1$, subcritical relative to the orbit
$\gamma^*$.} \label{figsubcritic2}
\end{figure}

To complete the analysis, let us consider orbits in $V(\alpha, 0)$, with $\alpha\in \left(\phi_1, \frac{\pi}{2}\right),$ 
as represented in Figure~\ref{figsubcritic2}.

%%%%%%%%%%%%%%%%%%%%%%%%%%%%%%%%%%%%%%%%%%%%%%%%%%%%%%%%%%%%%%%%%%%%%%%%
%
%
%
%
%
%
%
%
%
%
%
%
%%%%%%%%%%%%%%%%%%%%%%%%%%%%%%%%%%%%%%%%%%%%%%%%%%%%%%%%%%%%%%%%%%%%%%%%

It is clear from this plot that the time spent by this orbit is
$T_{C_{\ell}}(\alpha) := 2T(\alpha) - T_1(\alpha,\phi_0) + T_1(\alpha, \phi_1)$
(see subsection~\ref{313} for a justification of this notation $T_{C_{\ell}}$). 

It is also clear that $T_{C_{\ell}}(\alpha) \to T^*$ as $\alpha\to \phi_1.$ We shall prove that, for $\alpha>\phi_1$
sufficiently close to $\phi_1$, we have $T_{C_{\ell}}(\alpha) < T^*.$ This is not as easy to prove as in the previous cases.
We start by considering in \eqref{T} and \eqref{T1} a new variable $\tilde{\alpha} :=\sin^2\alpha,$ and changing in \eqref{T} the
integration variable $x\mapsto \theta$ where $\sin x = \sqrt{\tilde{\alpha}}\sin\theta.$ This allow us to write
\begin{equation}
 T_{C_{\ell}}(\alpha) =  \tilde{T}_{C_{\ell}}(\tilde{\alpha}) := \sqrt{2}\int_0^{\pi/2}\frac{d\theta}
 {\sqrt{1-\tilde{\alpha}\sin^2\theta}} +
 \frac{1}{\sqrt{2}}\int_{\phi_0}^{\phi_1}\frac{dx}{\sqrt{\tilde{\alpha}-\sin^2x}},  \label{time315}
\end{equation}
Differentiating we obtain
 \begin{eqnarray*}
  \frac{d\tilde{T}_{C_{\ell}}}{d\tilde{\alpha}} & = & 
  \frac{1}{\sqrt{2}}\int_0^{\pi/2}\frac{\sin^2\theta}{\left(1-\tilde{\alpha}\sin^2\theta\right)^{3/2}}d\theta -
  \frac{1}{2\sqrt{2}}\int_{\phi_0}^{\phi_1}\frac{dx}{\left(\tilde{\alpha}-\sin^2x\right)^{3/2}}, 
 \end{eqnarray*}
 and computing the second derivative we obtain
 \begin{eqnarray*}
  \frac{d^2\tilde{T}_{C_{\ell}}}{d\tilde{\alpha}^2} & = & 
  \frac{3}{2\sqrt{2}}\int_0^{\pi/2}\frac{\sin^4\theta}{\left(1-\tilde{\alpha}\sin^2\theta\right)^{5/2}}d\theta +
  \frac{3}{4\sqrt{2}}\int_{\phi_0}^{\phi_1}\frac{dx}{\left(\tilde{\alpha}-\sin^2x\right)^{5/2}}>0. 
 \end{eqnarray*}
Hence, $\tilde{T}_{C_{\ell}}$ is a convex function of $\tilde{\alpha} := \sin^2\alpha \in (\sin^2\phi_1, 1).$
From the definition of $T_{C_{\ell}}$ and $\tilde{T}_{C_{\ell}}$, the above expressions, and Proposition~\ref{monotonicity}, 
we also conclude that
$\tilde{T}_{C_{\ell}} \to +\infty$ as $\tilde{\alpha}\to 1,$ and $\frac{d\tilde{T}_{C_{\ell}}}{d\tilde{\alpha}} \to -\infty$
as $\tilde{\alpha}\to \sin^2\phi_1;$  however, note that $\tilde{T}_{C_{\ell}}\to T^*$ as $\tilde{\alpha}\to \sin^2\phi_1$
(see the start of this paragraph).

This behaviour obviously implies the existence of a single local extrema (a minimum) of $\tilde{T}_{C_{\ell}}$, 
and hence of $T_{C_{\ell}},$ in the interior of their respetive 
intervals of definition, and thus $T_{C_{\ell}}(\alpha) < T^*$ when $\alpha>\phi_1$
sufficiently close to $\phi_1.$ This justifies us calling this situation a (local) subcritical case relative to $\gamma^*$.
We emphasize that the situation is \emph{local}\/: if $\alpha$ is larger than the minimizer of $T_{C_{\ell}}(\alpha),$ the
value of this function increases without bound as $\alpha\to \pi/2,$ and thus,at some point, it will certainly be larger than $T^*$.

\begin{figure}[h]
\begin{center}
\psfrag{A}{$A$}
\psfrag{Cl}{$C_{\ell}$}
\psfrag{Cr}{$C_r$}
\psfrag{D}{$D$}
\psfrag{L}{$L$}
\psfrag{y-L}{$y(-L)$}
\psfrag{gdown}{$\gamma_*$}
\psfrag{gup}{$\gamma^*$}
\includegraphics[scale=0.30]{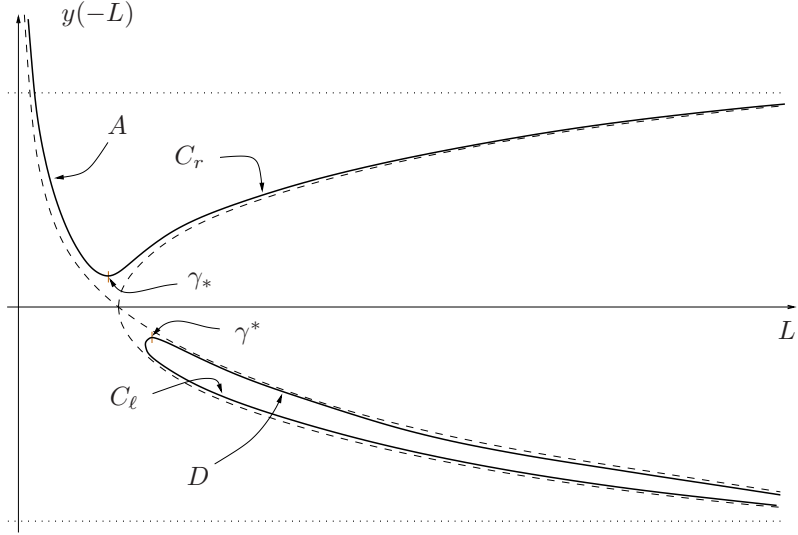}
\end{center}
\caption{Solid lines: portion of the bifurcation diagram when $\phi_0<\phi_1$ constructed from the analysis of the time-maps 
about $\gamma_*$ and $\gamma^*$ presented in subsections~\ref{311}--\ref{315}. 
Dashed lines: the corresponding diagram when $\phi_0=\phi_1$ (from \cite{CGGP}).
The designation of the orbits by letters $A$, $C_{\ell}$, $C_r$ and $D$ correspond to those used  
in \cite{CGGP}: see Table 1 and Figure 8 of that article.} \label{figbifdiagram1}
\end{figure}

Collecting the results obtained in the subsections~\ref{311}--\ref{315} we obtain the bifurcation diagram
in Figure~\ref{figbifdiagram1}. Remark that, due to the symmetry of the system, the value of $y(-L)$ of
the orbits $\gamma_*$ and $\gamma^*$ have the same absolute value (and different signs).

%%%%%%%%%%%%%%%%%%%%%%%%%%%%%%%%%%%%%%%%%%%%%%%%%%%%%%%%%%%%%%%%%%%%%%%%
%
%
%
%
%
%
%
%
%
%
%
%
%%%%%%%%%%%%%%%%%%%%%%%%%%%%%%%%%%%%%%%%%%%%%%%%%%%%%%%%%%%%%%%%%%%%%%%%

\subsubsection{Other solution branches} \label{316}

In addition to the solution branches studied above and represented in Fig.~\ref{figbifdiagram1},  
\eqref{eq1}--\eqref{eq2} has
an infinite number of solution families, each corresponding to orbits circling the origin a complete
$k$ number of times ($k=1, 2, \ldots).$ As in the cases studied above, it is convenient to start by considering orbits
corresponding to solutions that satisfy the additional boundary condition $y(L)=0,$ and, as before, we denote those 
orbits by a star, in this case by $\gamma_{*k}$ and $\gamma^*_k$. 
Although they do not correspond to bifurcating points, they are very useful in organizing or knowledge
about the solution branches. In Table~1 we present the orbit $\gamma_{*k}$ and those
which form a connected branch with it when $L$ changes from the value corresponding to $\gamma_{*k}$.
In Table~2 we present the analogous picture concerning the orbit $\gamma_k^{*}$.

\begin{center}
\begin{table}
\caption{Branch of solutions to (\ref{eq1})-(\ref{eq2}), with $\phi_0<\phi_1$, 
winding $k$ full times around $\mbox{\boldmath $0$}$ and containing the solution $\gamma_{*k}$ (For $k=0$ the
portion of the orbits with a thin trace should be disregarded.)}
\begin{minipage}{\textwidth}
\begin{tabular}{cc}
\hline\hline
 & \\ 
\multicolumn{1}{c}{ Orbit $\gamma_{\alpha,k}$} & 
\multicolumn{1}{c}{ Time taken by the orbit $\gamma_{\alpha,k}$}\\ 
% nova linha
\multicolumn{1}{c}{(winds $k$ times around $\mbox{\boldmath $0$}$)} &
\multicolumn{1}{c}{\mbox{}}\\ \hline\hline

\psfrag{lCr}{{\scriptsize{$A$}}}
\rule[-3mm]{0mm}{18mm}\mbox{\includegraphics[scale=0.45]{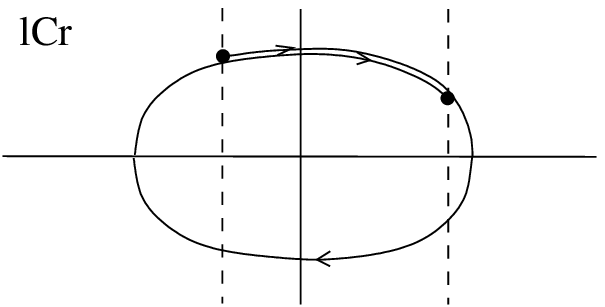}} & 
{ \rule[-3mm]{0mm}{2mm}$T_{A}(\alpha):=4kT(\alpha) + T_1(\alpha,\phi_0) + T_1(\alpha, \phi_1)$}\\
\hline
\psfrag{d}{{{$\gamma_{*k}$}}}
\rule[-3mm]{0mm}{18mm}\mbox{\includegraphics[scale=0.45]{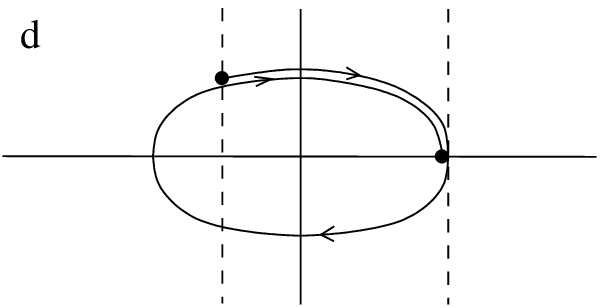}} & 
{ \rule[-3mm]{0mm}{2mm}$T_{*k}(\phi_1):=(4k+1)T(\phi_1) + T_1(\phi_1,\phi_0)$}\\
\hline
\psfrag{lCr}{{\scriptsize{$C_{r}$}}}
\rule[-3mm]{0mm}{18mm}\mbox{\includegraphics[scale=0.45]{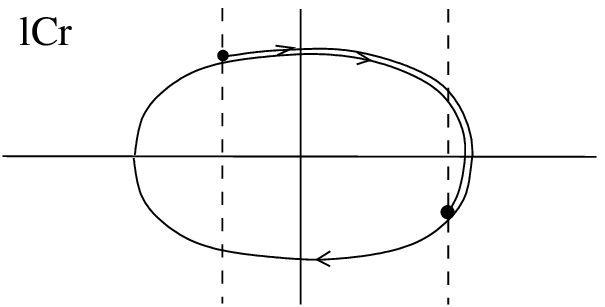}} & 
{ \rule[-3mm]{0mm}{2mm}$T_{C_{r}}(\alpha):=(4k+2)T(\alpha) + T_1(\alpha,\phi_0) -T_1(\alpha, \phi_1)$}\\
\hline\hline
\end{tabular}
\end{minipage}
\end{table}
\end{center}

\begin{center}
\begin{table}
\caption{Branch of solutions to (\ref{eq1})-(\ref{eq2}), with $\phi_0<\phi_1$, 
winding $k$ full times around $\mbox{\boldmath $0$}$ and containing the solution $\gamma_k^{*}$ (For $k=0$ the
portion of the orbits with a thin trace should be disregarded.)}
\begin{minipage}{\textwidth}
\begin{tabular}{cc}
\hline\hline
 & \\ 
\multicolumn{1}{c}{ Orbit $\gamma_{\alpha,k}$} & 
\multicolumn{1}{c}{ Time taken by the orbit $\gamma_{\alpha,k}$}\\ 
% nova linha
\multicolumn{1}{c}{(winds $k$ times around $\mbox{\boldmath $0$}$)} &
\multicolumn{1}{c}{\mbox{}}\\ \hline\hline

\psfrag{lCl}{{\scriptsize{$C_{\ell}$}}}
\rule[-3mm]{0mm}{18mm}\mbox{\includegraphics[scale=0.45]{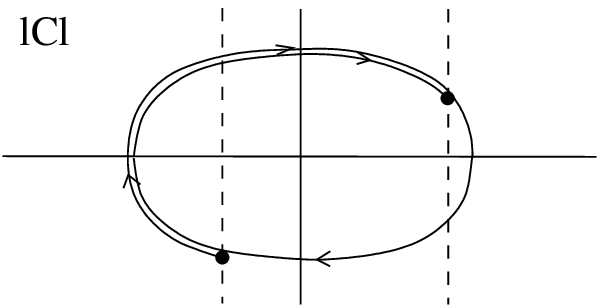}} & 
{ \rule[-3mm]{0mm}{2mm}$T_{C_{\ell}}(\alpha):=(4k+2)T(\alpha) - T_1(\alpha,\phi_0) + T_1(\alpha, \phi_1)$}\\
\hline
\psfrag{U}{{{$\gamma^*_{k}$}}}
\rule[-3mm]{0mm}{18mm}\mbox{\includegraphics[scale=0.45]{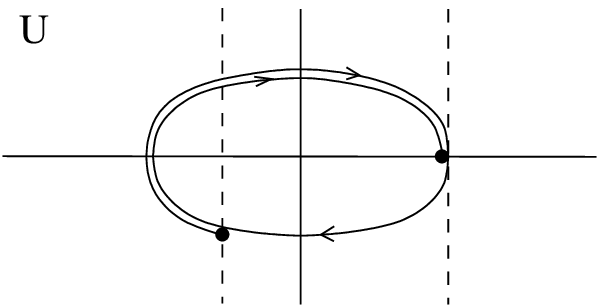}} & 
{ \rule[-3mm]{0mm}{2mm}$T_k^{*}(\phi_1):=(4k+3)T(\phi_1) - T_1(\phi_1,\phi_0)$}\\
\hline
\psfrag{D}{{\scriptsize{$D$}}}
\rule[-3mm]{0mm}{18mm}\mbox{\includegraphics[scale=0.45]{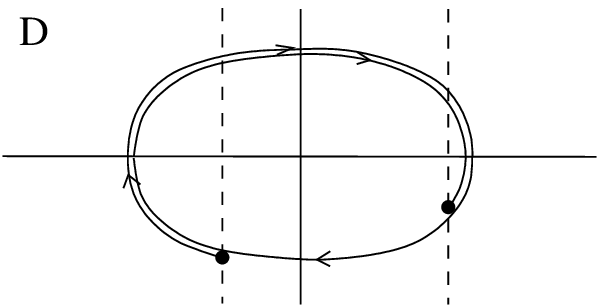}} & 
{ \rule[-3mm]{0mm}{2mm}$T_{D}(\alpha):=4(k+1)T(\alpha) - T_1(\alpha,\phi_0) -T_1(\alpha, \phi_1)$}\\
\hline\hline
\end{tabular}
\end{minipage}
\end{table}
\end{center}

Observe that these orbits are analogous to those studied in the previous subsections, which can be
considered the case $k=0$ in this description (i.e., the orbits do not complete a full turn around the origin).
The amounts of time spent by each of these orbits are exactly those of the corresponding ones in subsections~\ref{311}--\ref{315}
with the addition of $4kT(\alpha),$ which is the time of $k$ full turns about the origin.

The following conclusions are easily drawn:

% \begin{description}
% \item[a] From the definitions of the time-maps it follows that $T_{k}^{*}(\phi_1) < T_{(k+1)*}(\phi_1).$
% \end{description}
\mbox{} \\ 
\noindent
{\bf a:} From the definitions of the time-maps it follows that $T_{k}^{*}(\phi_1) < T_{(k+1)*}(\phi_1).$

\mbox{} \\ 
\noindent
{\bf b:} From \eqref{desigualdadeT} we immediately get $T_{k*}(\phi_1) < T_k^{*}(\phi_1).$

\mbox{} \\ 
\noindent
{\bf c:} From the results in subsections~\ref{313} and~\ref{314} and the fact that the time spent
by the orbits with $k>0$ is equal the to the time spent by those with $k=0$ plus $4kT(\alpha),$ we easily conclude
that $T_{C_r}(\alpha) > T_{*k}(\phi_1)$ and $T_D(\alpha)>T_k^{*}(\phi_1).$

\mbox{} \\ 
\noindent
{\bf d:} The study of the relation between $T_{C_{\ell}}(\alpha)$ and $T_k^*(\phi_1)$, for
$\alpha>\phi_1$ sufficiently close to $\phi_1$ proceeds exactly as in subsection~\ref{315}, paying attention to the
fact that we need to add $4kT(\alpha)$ to those computations. Since
$T'(\alpha)\to T(\phi_1) \in (0, +\infty)$ as $\alpha \downarrow \phi_1$, and
$\tilde{T}''(\tilde{\alpha})>0,$ the addition of $4kT(\alpha)$ to the right-hand side of \eqref{time315}
does nor change the conclusion. Hence we have $T_{C_{\ell}}(\alpha) < T_k^*(\phi_1),$ for $\alpha-\phi_1>0$
sufficiently small. Also, the other conclusions infered from the convexity of $\tilde{\alpha}\mapsto \tilde{T}_{C_{\ell}}
(\tilde{\alpha})$ remain valid.

\mbox{} \\ 
{\bf e:} Finally, it remains to study the relation between $T_{A}(\alpha)$ and $T_{*k}(\phi_1).$ The
analysis also follows that presented in subsection~\ref{315}. Changings variables as in subsection~\ref{315}
we can write an expression for $T_A(\alpha)$ similar to \eqref{time315}, namely
\begin{eqnarray*}
 \lefteqn{T_A(\alpha)  =  \tilde{T}_A(\tilde{\alpha}) := }  \\
 & := &  2\sqrt{2}k\int_0^{\pi/2}\frac{d\theta}{\sqrt{1-\tilde{\alpha}\sin^2\theta}}
 + \frac{1}{\sqrt{2}}\int_0^{\phi_0}\frac{dx}{\sqrt{\tilde{\alpha}-\sin^2x}}
 + \frac{1}{\sqrt{2}}\int_0^{\phi_0}\frac{dx}{\sqrt{\tilde{\alpha}-\sin^2x}}.
\end{eqnarray*}
Now, the convexity argument employed in subsection~\ref{315} and also used in case {\bf d} above, can again be applied to
conclude that, for $\alpha-\phi_1>0$ sufficiently small, type $A$ orbits satisfy $T_A(\alpha)< T_{*k}(\phi_1)$
and the corresponding branch in the diagram $L$ vs. $y(-L)$ is convex. Note that, in contrast to the case 
studied in subsection~\ref{315}, but as was the case in \cite{CGGP}, the branches of type $A$ solutions have 
a (unique, by convexity) saddle-node, since we know that, from Proposition~\ref{monotonicity}, $T_A(\alpha)\to +\infty$ as 
$\alpha\to \frac{\pi}{2}.$

%\mbox{} \\ 

\begin{figure}[h]
\begin{center}
\psfrag{A}{$A$}
\psfrag{Cl}{$C_{\ell}$}
\psfrag{Cr}{$C_r$}
\psfrag{D}{$D$}
\psfrag{L}{$L$}
\psfrag{y-L}{$y(-L)$}
\psfrag{gdown}{$\gamma_*$}
\psfrag{gup}{$\gamma^*$}
\psfrag{gdownk1}{$\gamma_{*1}$}
\psfrag{gupk1}{$\gamma^*_1$}
\psfrag{gdownk2}{$\gamma_{*2}$}
\psfrag{gupk2}{$\gamma^*_2$}
\includegraphics[scale=0.30]{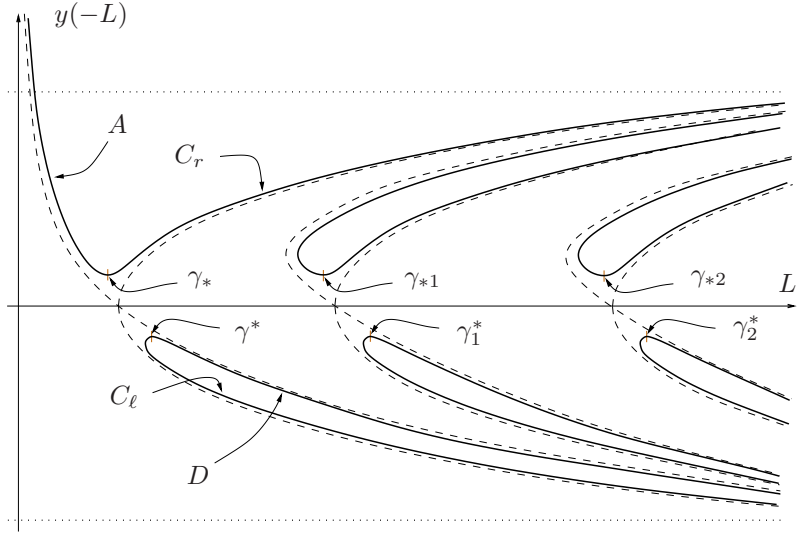}
\end{center}
\caption{Solid lines: portion of the bifurcation diagram when $\phi_0<\phi_1$ constructed from the analysis 
presented in subsections~\ref{311}--\ref{316}. 
Dashed lines: the corresponding diagram when $\phi_0=\phi_1$ (from \cite{CGGP}).
The designation of the orbits by letters $A$, $C_{\ell}$, $C_r$ and $D$ correspond to those used  
in \cite{CGGP}: see Table 1 and Figure 8 of that article.} \label{figbifdiagram2}
\end{figure}

Thus, we conclude from these results that, for each $k$, the relation of the various types of orbits among themselves is
the same as existed in the case $k=0$ illustrated in Figure~\ref{figbifdiagram1}. 
We collect the results obtained thus far in the bifurcation diagram of Figure~\ref{figbifdiagram2}.
Observe that, due to the symmetry of the system, the value of $y(-L)$ of
the orbits $\gamma_{*k}$ are the same for all $k$, and the same happens for $\gamma^*_k;$
as was the case when $k=0$, for all $k$ these values in $\gamma_{*k}$ and in $\gamma^*_k$  
have the same absolute value (and different signs).

%%%%%%%%%%%%%%%%%%%%%%%%%%%%%%%%%%%%%%%%%%%%%%%%%%%%%%%%%%%%%%%%%%%%%%%%
%
%
%
%
%
%
%
%
%
%
%
%
%%%%%%%%%%%%%%%%%%%%%%%%%%%%%%%%%%%%%%%%%%%%%%%%%%%%%%%%%%%%%%%%%%%%%%%%

\subsection{Case $\phi_0>\phi_1$}\label{32}

\begin{figure}[h]
\begin{center}
\psfrag{a}{\scriptsize{$(a)$}}
\psfrag{b}{\scriptsize{$(b)$}}
\psfrag{x}{$x$}
\psfrag{y}{$y$}
\psfrag{-p}{\scriptsize{$x=-\phi_0$}}
\psfrag{p}{\scriptsize{$x=\phi_1$}}
\psfrag{-pi}{$-\frac{\pi}{2}$}
\psfrag{pi}{$\frac{\pi}{2}$}
\psfrag{ga}{$\gamma_*$}
\psfrag{gb}{$\gamma^*$}
\includegraphics[scale=0.23]{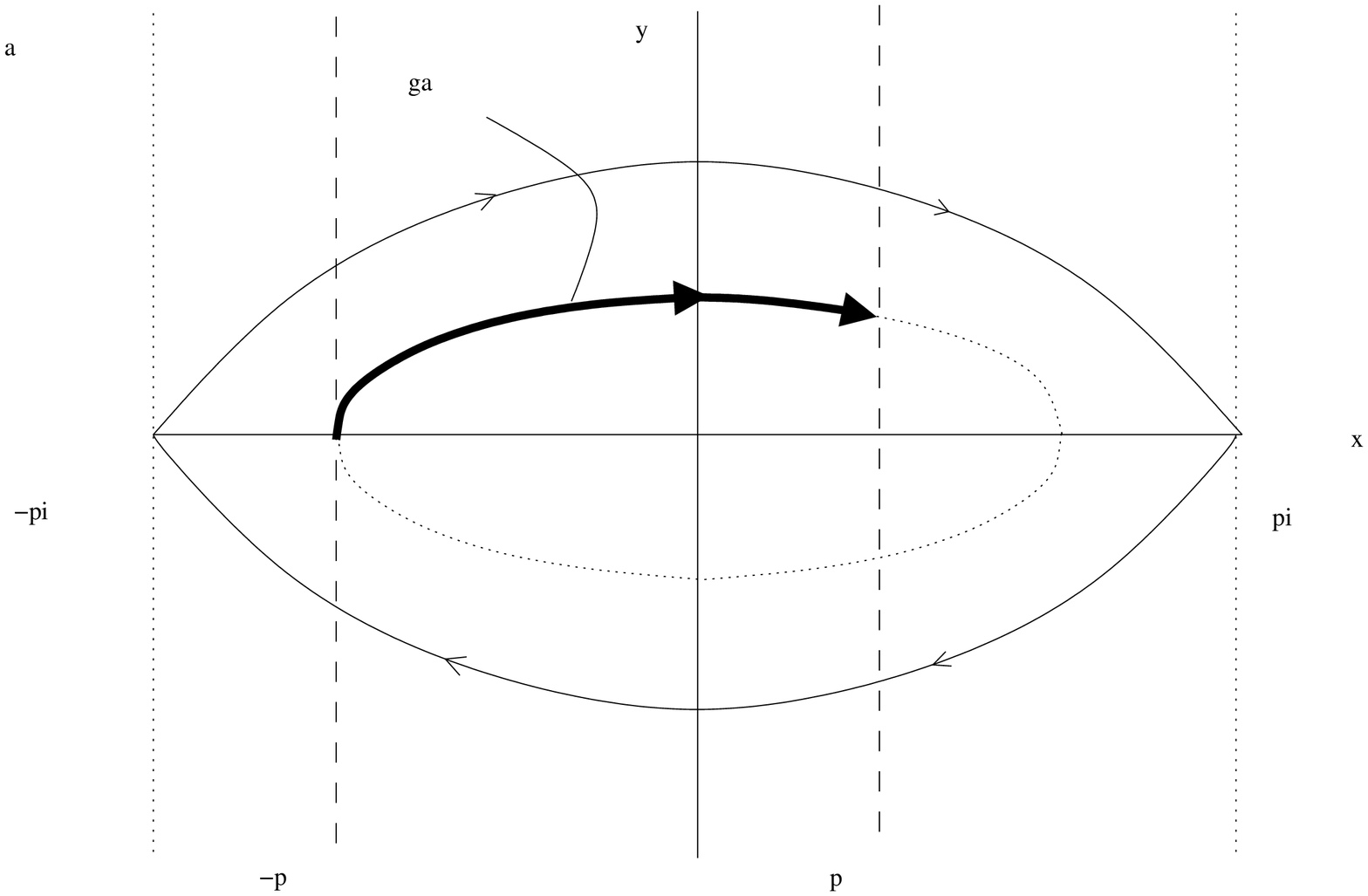}\qquad\includegraphics[scale=0.23]{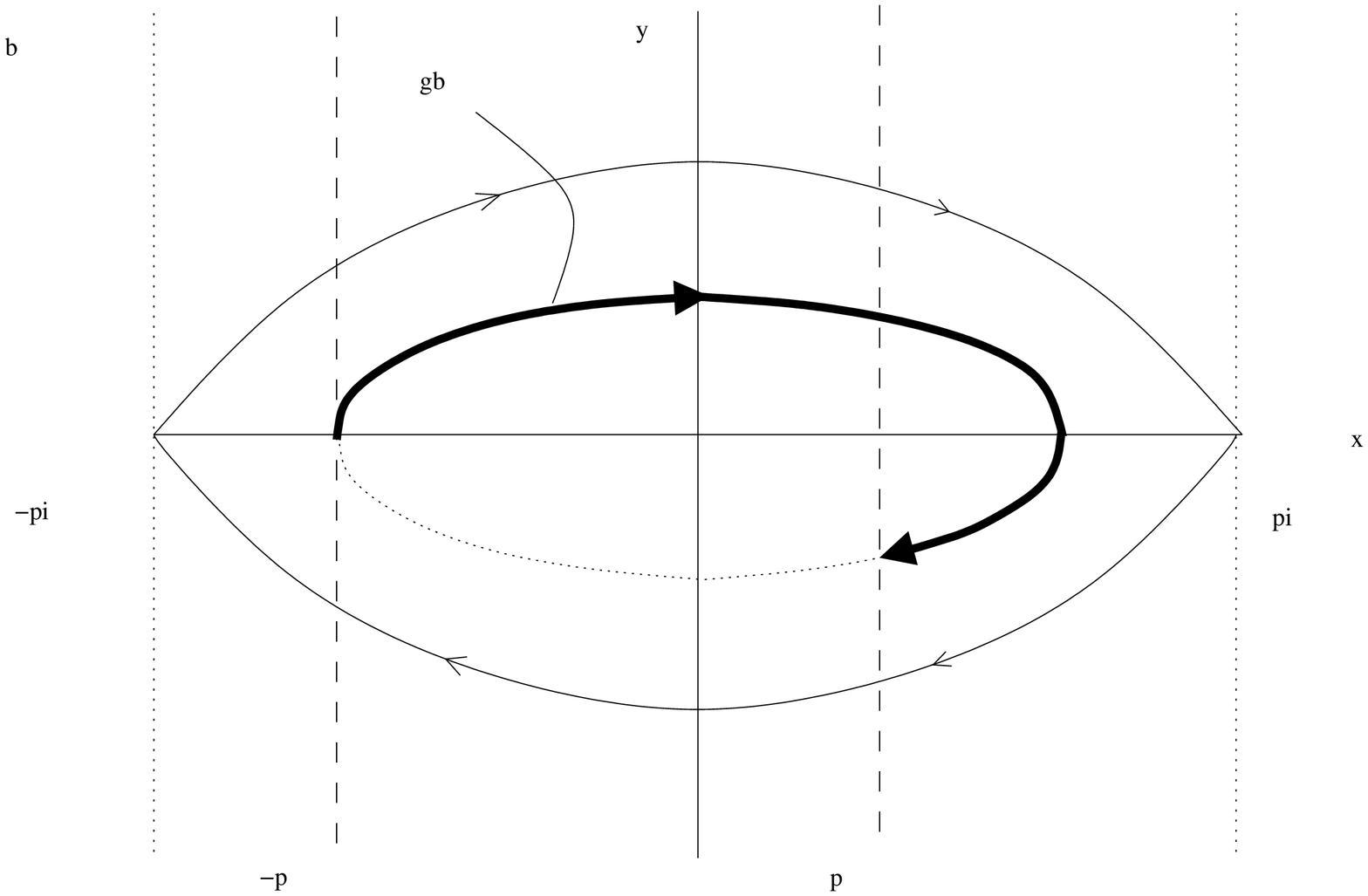}
\end{center}
\caption{Two ``critical'' orbits of (\ref{eq1})-(\ref{eq2}) with $\phi_0>\phi_1$: (a) $\gamma_*$ 
 when $2L=T_*:=T(\phi_1)+T_1(\phi_1,\phi_0)$;
 (b) $\gamma^*$ 
 when $2L=T^*:=3T(\phi_1)-T_1(\phi_1,\phi_0)$} \label{figcritic2}
\end{figure}

The analysis of the case $\phi_0>\phi_1$ proceeds in a way entirely similar to the case $\phi_0<\phi_1$
and so we will not present the details of the arguments in what follows. We will, in the next figures, exhibit the
plots of the several types of orbits and the bifurcation diagram obtained. We start, in Figure~\ref{figcritic2}, 
by the orbits that, at $t=-L$,
satisfy the additional boundary condition $y(-L)=0$, which we designate by ``critical'' orbits, as done in the similar
situation in subsection~\ref{311}.

%%%%%%%%%%%%%%%%%%%%%%%%%%%%%%%%%%%%%%%%%%%%%%%%%%%%%%%%%%%%%%%%%%%%%%%%
%
%
%
%
%
%
%
%
%
%
%
%
%%%%%%%%%%%%%%%%%%%%%%%%%%%%%%%%%%%%%%%%%%%%%%%%%%%%%%%%%%%%%%%%%%%%%%%%

\begin{figure}[h]
\begin{center}
\psfrag{a}{\scriptsize{$(b)$}}
\psfrag{a2}{\scriptsize{$(a)$}}
\psfrag{b}{\scriptsize{$(b)$}}
\psfrag{x}{$x$}
\psfrag{y}{$y$}
\psfrag{-p}{\scriptsize{$x=-\phi_0$}}
\psfrag{p}{\scriptsize{$x=\phi_1$}}
\psfrag{-pi}{$-\frac{\pi}{2}$}
\psfrag{pi}{$\frac{\pi}{2}$}
\psfrag{gb}{$\gamma_*$}
\includegraphics[scale=0.23]{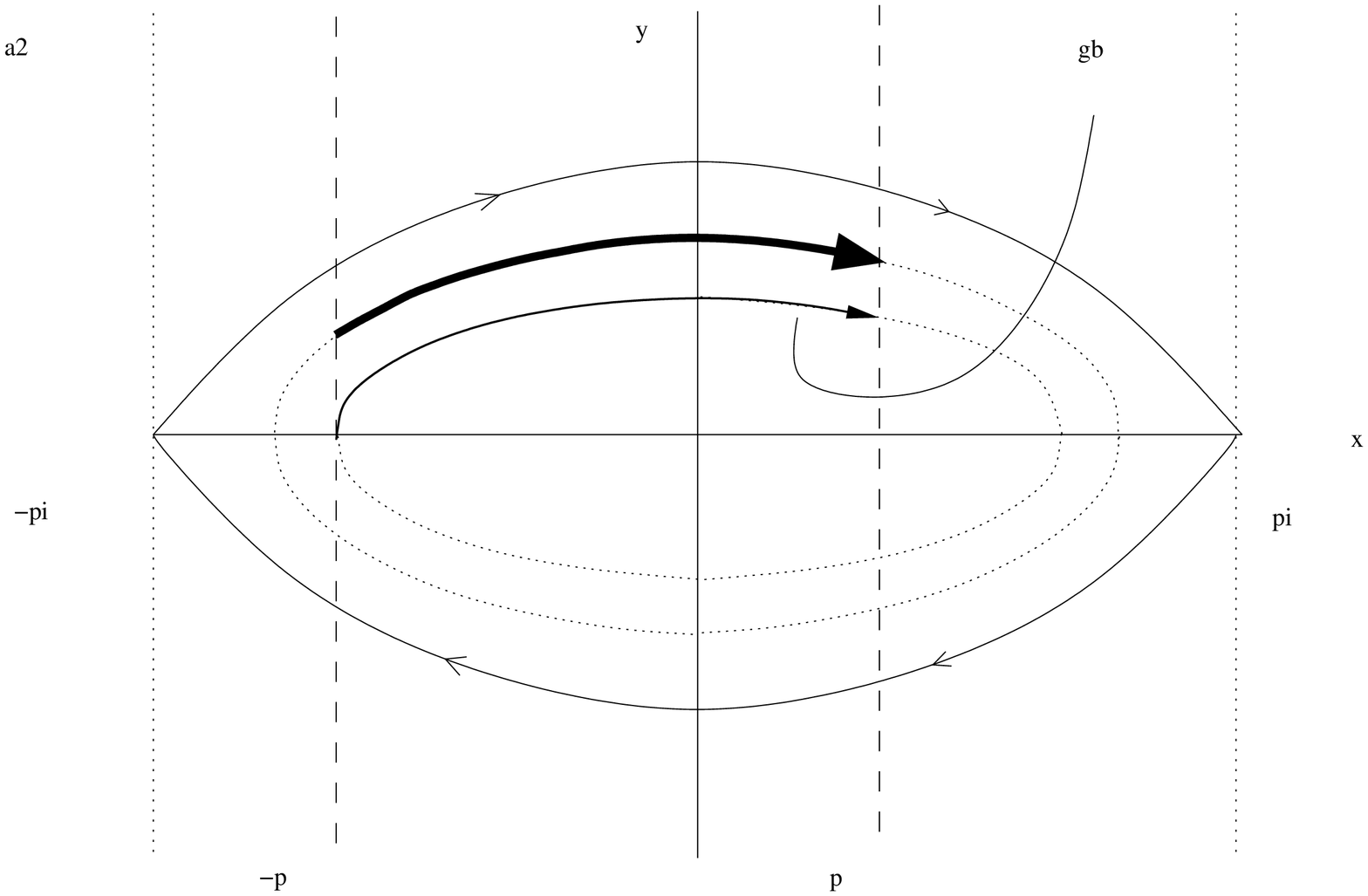}\qquad\includegraphics[scale=0.23]{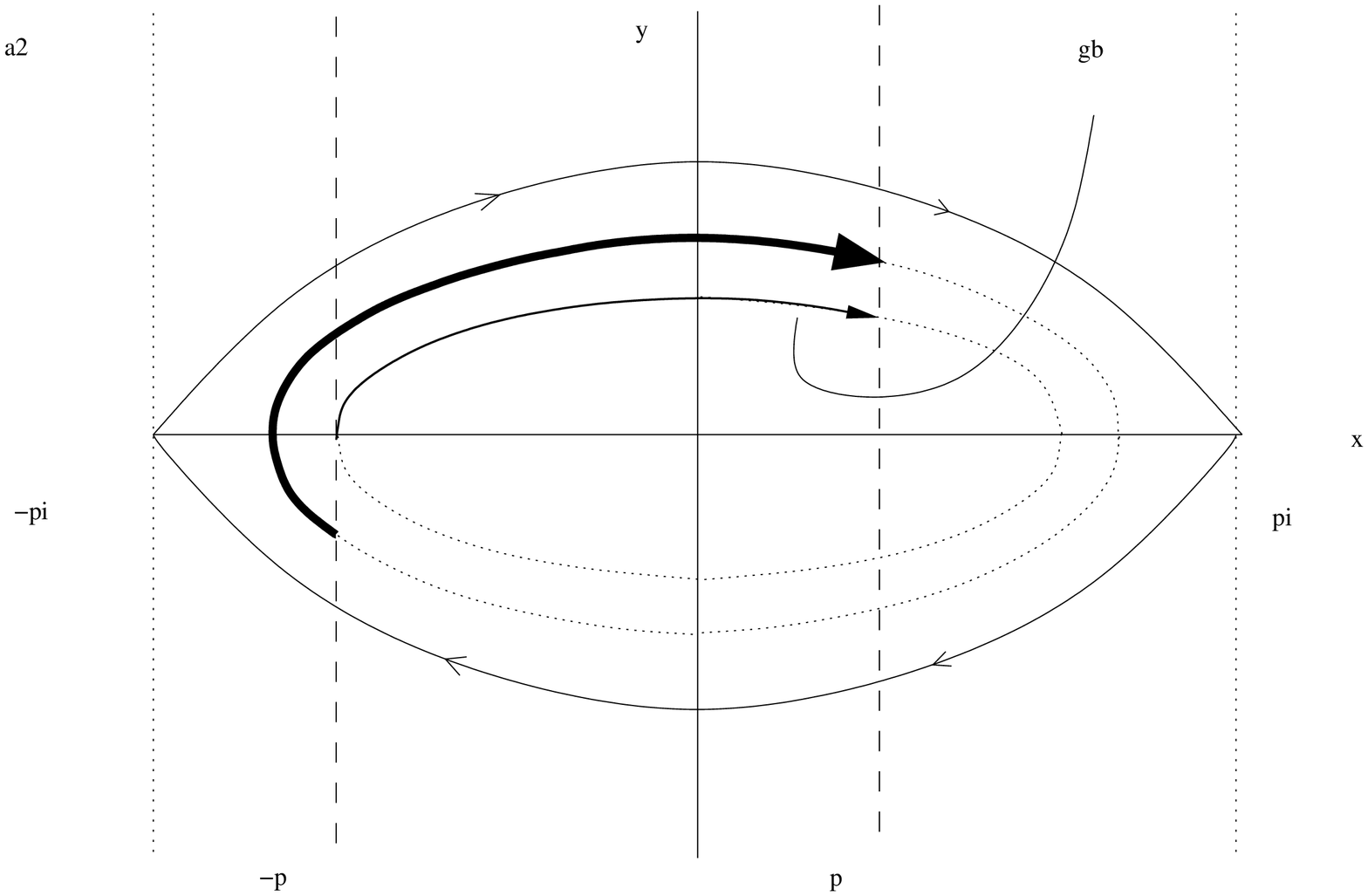}
\end{center}
\caption{Orbits of (\ref{eq1})-(\ref{eq2}), with $\phi_0>\phi_1$ which are: (a) subcritical relative to the orbit $\gamma_*$; 
(b) supercritical relative to the orbit $\gamma_*$.} \label{fig_branch_lower_*}
\end{figure}

\begin{figure}[h]
\begin{center}
\psfrag{a}{\scriptsize{$(b)$}}
\psfrag{a2}{\scriptsize{$(a)$}}
\psfrag{b}{\scriptsize{$(b)$}}
\psfrag{x}{$x$}
\psfrag{y}{$y$}
\psfrag{-p}{\scriptsize{$x=-\phi_0$}}
\psfrag{p}{\scriptsize{$x=\phi_1$}}
\psfrag{-pi}{$-\frac{\pi}{2}$}
\psfrag{pi}{$\frac{\pi}{2}$}
\psfrag{ga}{$\gamma_*$}
\psfrag{gb}{$\gamma^*$}
\includegraphics[scale=0.23]{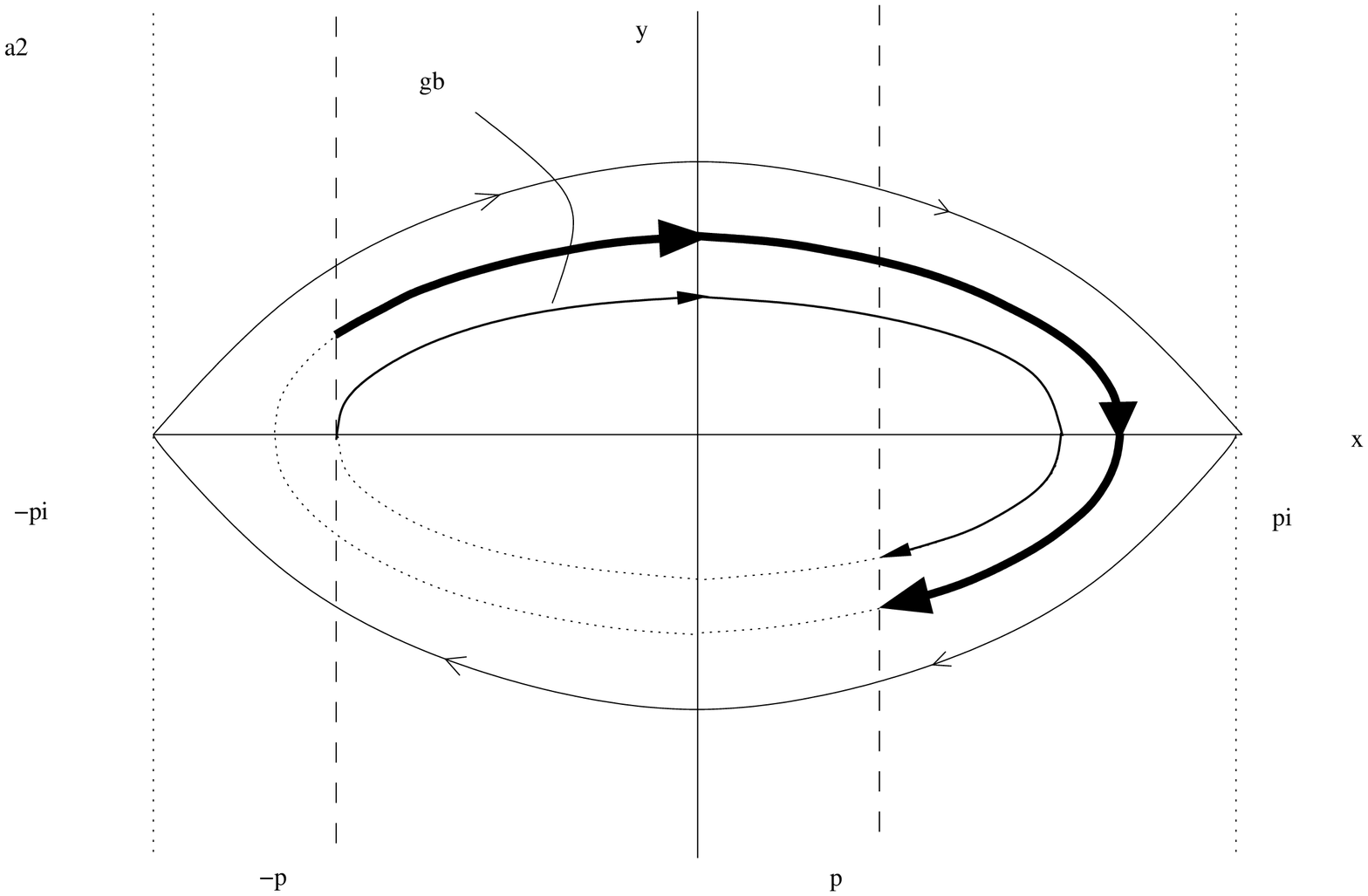}\qquad\includegraphics[scale=0.23]{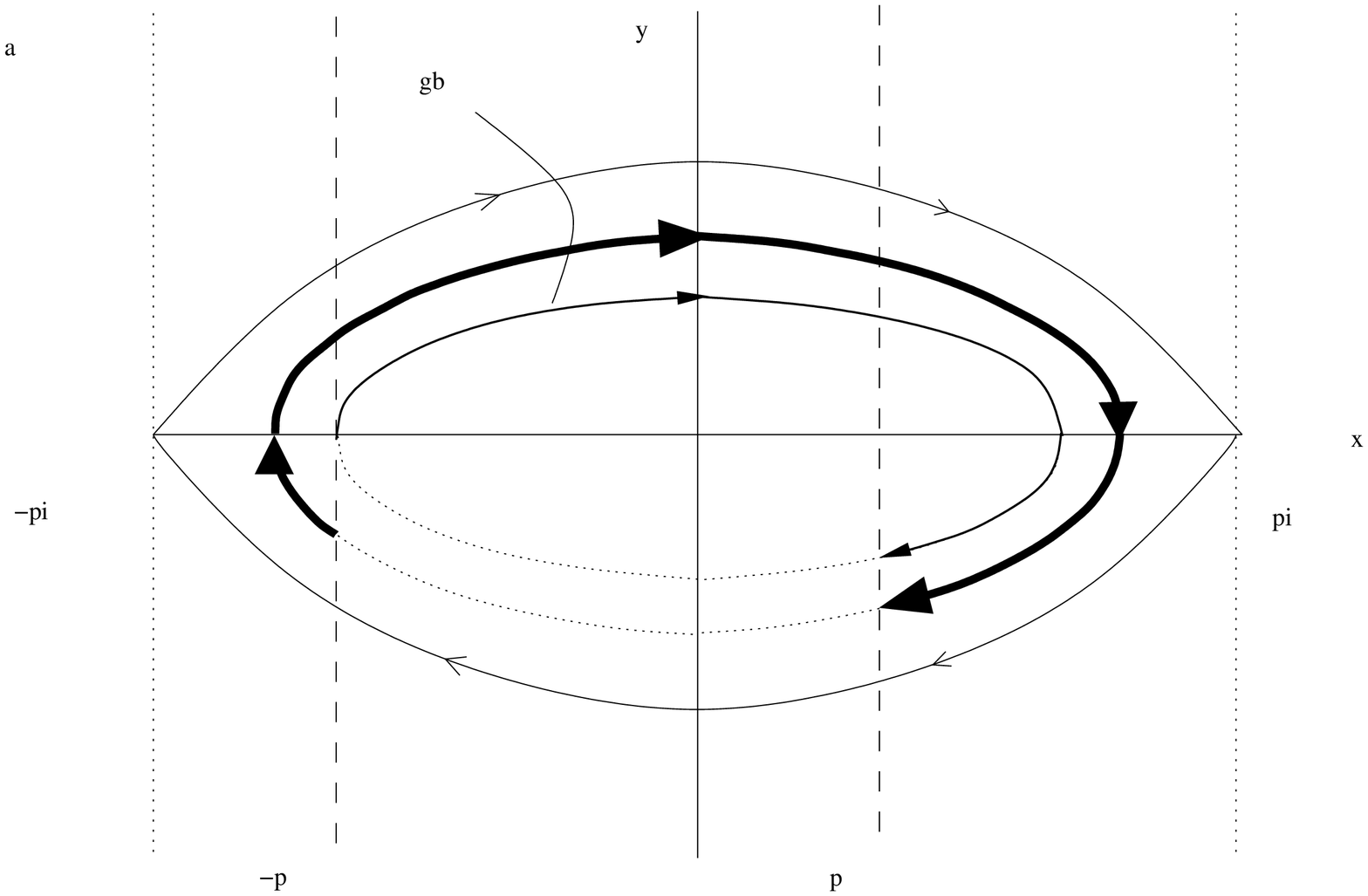}
\end{center}
\caption{Orbits of (\ref{eq1})-(\ref{eq2}), with $\phi_0>\phi_1$ which are: (a) subcritical relative to the orbit $\gamma^*$ (for
$\alpha-\phi_0>0$ sufficiently small); 
(b) supercritical relative to the orbit $\gamma^*$.} \label{fig_branch_upper_*}
\end{figure}

Due to the symmetry of the problem relative to the transformations $x\mapsto -x$ and $\phi_0 \leftrightarrow \phi_1$,
we conclude that, from each ``critical'' orbit emerges two branches, a subcritical and a supercritical, with exactly the same 
properties as obtained for the corresponding branches in subsections~\ref{312}--\ref{315}.  
These orbits are illustrated in figures~\ref{fig_branch_lower_*} and~\ref{fig_branch_upper_*}.

%%%%%%%%%%%%%%%%%%%%%%%%%%%%%%%%%%%%%%%%%%%%%%%%%%%%%%%%%%%%%%%%%%%%%%%%
%
%
%
%
%
%
%
%
%
%
%
%
%%%%%%%%%%%%%%%%%%%%%%%%%%%%%%%%%%%%%%%%%%%%%%%%%%%%%%%%%%%%%%%%%%%%%%%%

%%%%%%%%%%%%%%%%%%%%%%%%%%%%%%%%%%%%%%%%%%%%%%%%%%%%%%%%%%%%%%%%%%%%%%%%
%
%
%
%
%
%
%
%
%
%
%
%
%%%%%%%%%%%%%%%%%%%%%%%%%%%%%%%%%%%%%%%%%%%%%%%%%%%%%%%%%%%%%%%%%%%%%%%%

In an entirely analogous way to what was presented in section~\ref{316}, we also have the solution branches 
corresponding to orbits circling the origin a complete number $k\geqs 1$ of turns. 

\begin{figure}[h]
\begin{center}
\psfrag{A}{$A$}
\psfrag{Cl}{$C_{\ell}$}
\psfrag{Cr}{$C_r$}
\psfrag{D}{$D$}
\psfrag{L}{$L$}
\psfrag{y-L}{$y(-L)$}
\psfrag{gdown}{$\gamma_*$}
\psfrag{gup}{$\gamma^*$}
\psfrag{gdownk1}{$\gamma_{*1}$}
\psfrag{gupk1}{$\gamma^*_1$}
\psfrag{gdownk2}{$\gamma_{*2}$}
\psfrag{gupk2}{$\gamma^*_2$}
\includegraphics[scale=0.30]{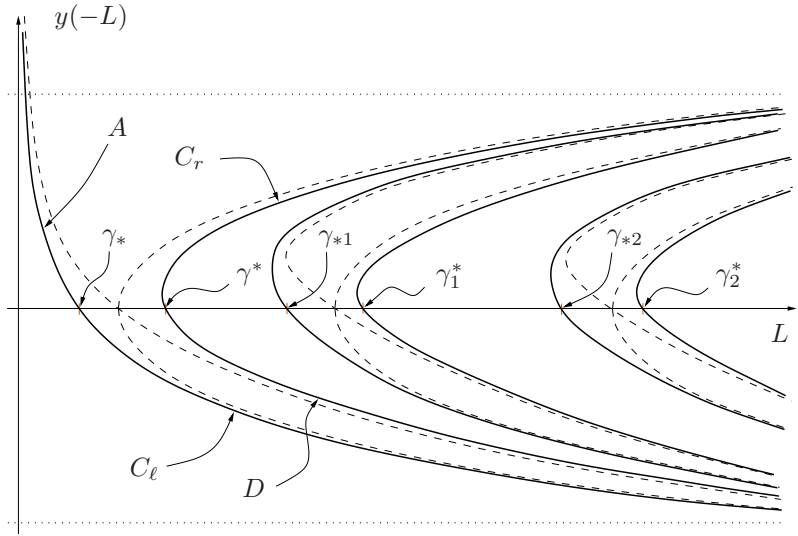}
\end{center}
\caption{Solid lines: portion of the bifurcation diagram when $\phi_0>\phi_1$ constructed from what was presented 
and discussed in subsection~\ref{32}. 
Dashed lines: the corresponding diagram when $\phi_0=\phi_1$ (from \cite{CGGP}).
The designation of the orbits by letters $A$, $C_{\ell}$, $C_r$ and $D$ correspond to those used  
in \cite{CGGP}: see Table 1 and Figure 8 of that article.} \label{figbifdiagram3}
\end{figure}

Collecting these results we can plot the bifurcation diagram corresponding to the case $\phi_0>\phi_1$. This is done
in Figure~\ref{figbifdiagram3}. To understand the apparently drastic difference relative to the diagram
for the case $\phi_0<\phi_1$ presented in Figure~\ref{figbifdiagram2} we need to bear in mind the fact that in \emph{both} cases
what is being plot in the vertical axis is the value of $y(-L)$ of the corresponding orbit. If, in the case $\phi_0>\phi_1$,
we choose to plot the value of $y(L)$ instead, by the symmetry considerations alluded to above, the corresponding bifurcation diagram 
will be equal to that of Figure~\ref{figbifdiagram2}.

%%%%%%%%%%%%%%%%%%%%%%%%%%%%%%%%%%%%%%%%%%%%%%%%%%%%%%%%%%%%%%%%%%%%%%%%
%
%
%
%
%
%
%
%
%
%
%
%
%%%%%%%%%%%%%%%%%%%%%%%%%%%%%%%%%%%%%%%%%%%%%%%%%%%%%%%%%%%%%%%%%%%%%%%%

\section{Stability analysis of the equilibria}

In this section we present a brief study of the stability of the equilibria using the approach
of Maginu \cite{maginu78}. 
We believe it is possible,
by a modification of methods originally developed for homogeneous Neumann boundary conditions (see, e.g., \cite{fr}
and \cite[Section 4.3]{hale}) to provide more detailed information about the \emph{unstable}\/ solutions, in particular 
clarifying, for each unstable equilibrium, which directions are unstable, and to characterize 
their heteroclinic connections. This will be postponed to a latter work.

We pretend to classify as stable, asymptotically stable, or unstable the branches of equilibria determined in the last section.
The results of \cite{maginu78} relevant to our case are the theorems 3.1, 3.2 and 3.3. 
What the first two of these theorems state is that solutions $(x(t),y(t))$ of (\ref{eq1})--(\ref{eq2})
are asymptotically stable (as stationary solutions of the corresponding 
partial differential equation (\ref{eq1+})--(\ref{eq1##})) if $y(t)$ has no zeros in $[-L, L)$ or in $(-L,L]$; and it
is unstable if $y(t)$ has two or more zeros. 

Clearly, these results take care of the stability characterization of all the branches of solutions with $k\geqs 1$ (they are all unstable),
and also when $k=0$ of the branch denoted by $A$ (which is asymptotically stable), and by $D$ (which is unstable).

Theorem 3.3 of \cite{maginu78} is one of a series or results characterizing the case when $y(t)$ has a single zero in $[-L, L],$
 located in $(-L, L)$.  Maginus' result states that such an equilibrium $E$
is asymptotically stable if the corresponding time map $T_E(\alpha)$ is strictly increasing, and is unstable if it is 
strictly decreasing.

Applied to our case, this result will allow us to determine the stability of the
remaining cases, namely:  the branches $C_r$ and $C_{\ell}$ when $k=0$. 

Consider $\phi_0<\phi_1$. Let us start with the $C_r$ branch.
Clearly such solutions are of the type considered in \cite[Theorem 3.3]{maginu78} (the existence of 
a single time instant for which $y(t)=0$). In Subsection~\ref{313} we concluded that $\frac{dT_{C_r}}{d\alpha} >0.$
Hence, Maginu's result imply the branch is asymptotically stable.

Let us consider now the case of the $C_\ell$ branch. The relevant computations are the ones in subsection~\ref{315}, 
where we concluded that $T_{C_\ell}(\alpha)$ is convex, with a single local minimum. This means that 
$\frac{dT_{C_\ell}}{d\alpha} <0$ for the part of the $C_\ell$ branch to the left of the $\gamma^*$ and to the right
of the leftmost point of the branch, i.e., the saddle-node bifurcation point (which corresponds 
to the orbit for which  $T_{C_\ell}(\alpha)$
attains its unique minimum.) So, by \cite[Theorem 3.3]{maginu78}, these equilibria are unstable. For the
remaining part of the $C_\ell$ branch, i.e, 
for points of the orbit below the saddle-node bifurcation point, we have $\frac{dT_{C_\ell}}{d\alpha} >0$, and thus, again 
by \cite[Theorem 3.3]{maginu78}, the corresponding equilibria are asymptotically stable.

\begin{figure}[h]
\begin{center}
\psfrag{A}{$A$}
\psfrag{Cl}{$C_{\ell}$}
\psfrag{Cr}{$C_r$}
\psfrag{D}{$D$}
\psfrag{L}{$L$}
\psfrag{sn}{$SN$}
\psfrag{s}{$s$}
\psfrag{u}{$u$}
\psfrag{gdown}{$\gamma_*$}
\psfrag{gup}{$\gamma^*$}
\includegraphics[scale=0.76]{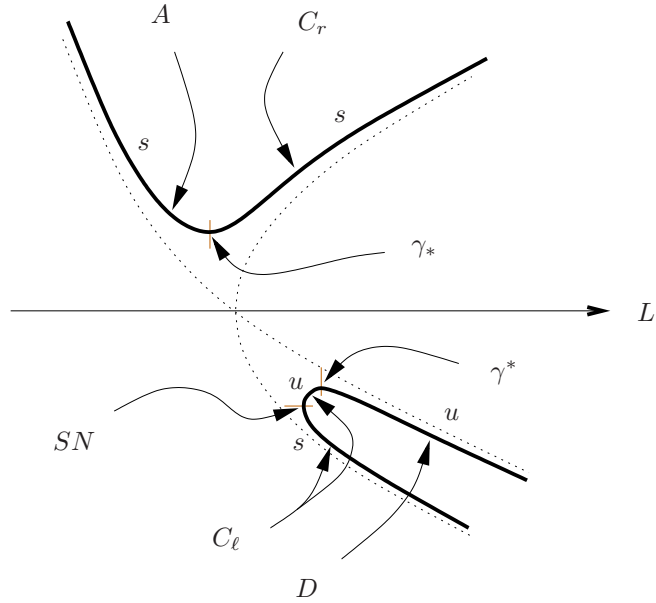}
\end{center}
\caption{Enlargement of the bifurcation diagram of Figure~\ref{figbifdiagram1}. 
The saddle-node bifurcation point refered to in the text is denoted by $SN$, and
the letters $s$ and $u$ denote branches of stable and unstable solutions, respectively.
The remaining notation is as in Figure~\ref{figbifdiagram1}.} \label{figbifdiagram1_exp}
\end{figure}

These stability conclusions for the $k=0$ branches are collected in Figure~\ref{figbifdiagram1_exp}.

%%%%%%%%%%%%%%%%%%%%%%%%%%%%%%%%%%%%%%%%%%%%%%%%%%%%%%%%%%%%%%%%%%%%%%%%
%
%
%
%
%
%
%
%
%
%
%
%
%%%%%%%%%%%%%%%%%%%%%%%%%%%%%%%%%%%%%%%%%%%%%%%%%%%%%%%%%%%%%%%%%%%%%%%%

Exactly the same results can be applied to the case when $\phi_0>\phi_1$ with analogous results: by
theorems 3.1 and 3.2 of \cite{maginu78} all the $k\geqs 1$
branches are unstable, as well as the $D$ branch, whereas the $A$ branch is asymptotically stable. An 
analysis corresponding to that in subsections~\ref{313} and~\ref{315} and the
application of theorem 3.3 of \cite{maginu78} results in the conclusion that 
$C_\ell$ is an asymptotically stable branch, and the portion of the $C_r$ branch between $y^*$ and the
leftmost point (a saddle-node) of the branch corresponds to unstable equilibria, whereas the points above this last point are
asymptotically stable equilibria.

\begin{figure}[h]
\begin{center}
\psfrag{A}{$A$}
\psfrag{Cl}{$C_{\ell}$}
\psfrag{Cr}{$C_r$}
\psfrag{D}{$D$}
\psfrag{L}{$L$}
\psfrag{sn}{$SN$}
\psfrag{s}{$s$}
\psfrag{u}{$u$}
\psfrag{gdown}{$\gamma_*$}
\psfrag{gup}{$\gamma^*$}
\includegraphics[scale=0.56]{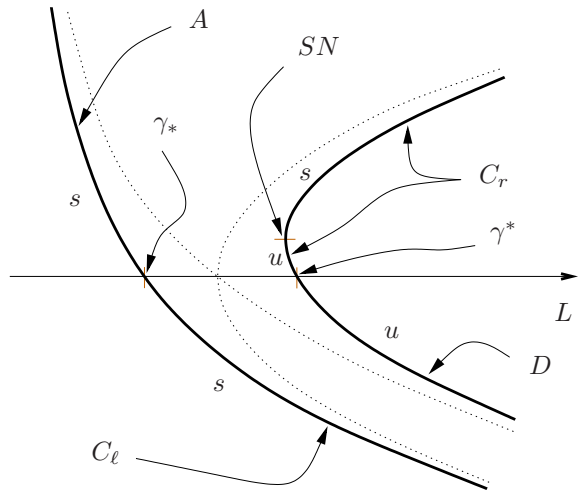}
\end{center}
\caption{Enlargement of the bifurcation diagram of Figure~\ref{figbifdiagram3}. 
The saddle-node bifurcation point refered to in the text is denoted by $SN$, and
the letters $s$ and $u$ denote branches of stable and unstable solutions, respectively.
The remaining notation is as in Figure~\ref{figbifdiagram3}.} \label{figbifdiagram3_exp}
\end{figure}

These  conclusions about the stability of the $k=0$ branches are collected in Figure~\ref{figbifdiagram3_exp}.

%%%%%%%%%%%%%%%%%%%%%%%%%%%%%%%%%%%%%%%%%%%%%%%%%%%%%%%%%%%%%%%%%%%%%%%%
%
%
%
%
%
%
%
%
%
%
%
%
%%%%%%%%%%%%%%%%%%%%%%%%%%%%%%%%%%%%%%%%%%%%%%%%%%%%%%%%%%%%%%%%%%%%%%%%

% \begin{acknowledgments} 
% F.P. da C. and J.T.P. were partially supported by FCT/Portugal through
% UID/MAT/04459/2013.
% Portions of Section 3.1 of this paper first appeared in the dissertation submited by M.I.M 
%    to the \emph{Universidad Nacional de Educaci\'on a Distancia,\/} 
%   Madrid, Spain, in June 2014, as part of the requirements for the  degree of
%   \emph{M\'aster en Matem\'aticas Avanzadas.\/}
% \end{acknowledgments}

%===================================================================================
%                                          References
%===================================================================================

\bibliographystyle{amsplain}

\end{document}